\documentclass[10pt,a4paper]{article}
\usepackage{amsmath,amsthm,amssymb,latexsym,epic}

\newtheorem{theorem}{Theorem}
\newtheorem{remark}[theorem]{Remark}
\newtheorem{lemma}[theorem]{Lemma}
\newtheorem{corollary}[theorem]{Corollary}
\newtheorem{proposition}[theorem]{Proposition}

\usepackage[all]{xy}
\usepackage{enumerate}
\newcommand{\IP}{\mathcal{I}^{\ast}}
\newcommand{\rank}{\mathrm{rank}}
\newcommand{\Sym}{\mathcal{S}}
\newcommand{\GR}{\mathcal{R}}
\newcommand{\GL}{\mathcal{L}}
\newcommand{\GH}{\mathcal{H}}
\newcommand{\GD}{\mathcal{D}}
\newcommand{\GJ}{\mathcal{J}}

\newcommand{\oIS}{\widetilde{\mathcal{I}}}
\newcommand{\oS}{\widetilde{\mathcal{S}}}
\newcommand{\oIP}{\widetilde{\mathcal{I}^{\ast}}}

\newcommand{\codom}{\mathrm{codom}}
\newcommand{\coran}{\mathrm{coran}}
\newcommand{\corank}{\mathrm{corank}}
\newcommand{\St}{\mathrm{St}}

\newcommand{\IS}{\mathcal{I}}

\newcommand{\Aut}{\mathrm{Aut}}

\newcommand{\dom}{\mathrm{dom}}
\newcommand{\ran}{\mathrm{ran}}

\newcommand{\lmod}{\mid\!}
\newcommand{\rmod}{\!\mid}

\newcommand{\N}{\mathcal{N}}
\newcommand{\PIP}{\mathcal{PI}^{\ast}}
\newcommand{\wPIP}{\overline{\mathcal{PI}^{\ast}}}
\renewcommand{\S}{\mathcal{S}}
\newcommand{\C}{\mathcal{C}}
\newcommand{\B}{\mathcal{B}}
\newcommand{\PB}{\mathcal{PB}}

\title{Two generalisations of the symmetric inverse semigroups}
\author{Ganna Kudryavtseva and Victor Maltcev}
\date{}
\begin{document}

\maketitle

\begin{abstract}
We introduce two generalisations of the full symmetric inverse 
semigroup $\IS_X$ and its dual semigroup $\IP_X$ -- inverse
semigroups $\PIP_X$ and $\wPIP_X$. Both of them have the same
carrier and contain $\IS_X$. Binary operations on $\PIP_X$ and
$\wPIP_X$ are reminiscent of the multiplication in $\IS_X$. We use a
convenient geometric way to realise elements from these two
semigroups. This enables us to study efficiently their inner
properties and to compare them with the corresponding properties of
$\IS_X$ and $\IP_X$.
\end{abstract}


\section{Introduction} \label{sec:intro}

One of the most natural examples of proper inverse semigroups (i.e.,
except groups) is the symmetric inverse semigroup $\IS_X$. Beside
pure combinatorial interest in this semigroup, it plays an important
role for the class of all inverse semigroups similar to that played
by the symmetric group $\S_X$ for the class of all groups. For some
facts about semigroup and combinatorial properties of $\IS_X$ we
refer the reader to~\cite{GM1}.

Seeking for further natural examples of inverse semigroups,
FitzGerald and Leech~\cite{FL}, using categorical methods, introduced
the \emph{dual symmetric inverse semigroup} $\IP_X$. Using more
general categorical approach, $\IP_X$ also appeared in~\cite{KM2}.
This semigroup also has a useful geometric realisation, which was
exploited in~\cite{EEF,Mal} to study some inner properties of
$\IP_X$.

In a recent work~\cite{KM1} there was found a new, representation
theoretic, link between $\IS_X$ and $\IP_X$.

In addition, both $\IP_X$ and $\IS_X$ belong to the class of the
so-called \emph{partition semigroups} \cite{Maz,Xi} and are
contained in the ``biggest partition semigroup'' $\C_X$ (see
Section~\ref{sec:def} for details). The latter semigroup was studied
mainly in the context of representation theory and cellular
algebras~\cite{HR, Jo, Ma, Xi}. Some pure semigroup aspects of
$\C_X$ were studied in~\cite{HR, Maz}.

In the present paper we aim at constructing two inverse semigroups
$\PIP_X$ and $\wPIP_X$, which are strongly related to $\IS_X$ and
$\IP_X$, though have much more complicated structure. We give transparent geometric definitions for these two
semigroups and then study their inner properties, focusing on
combinatorial aspects and their resemblance to $\IS_X$ and $\IP_X$.

The semigroups $\PIP_X$ and $\wPIP_X$ are natural also from the
representation theoretic point of view: $\wPIP_X$ is contained in a
bigger semigroup, the ``deformation'' of $\C_X$, whose semigroup
algebra naturally arises in the representation theory, see, e.g.,~\cite{HR}.
Some other representation theoretic aspects, where $\PIP_X$ and
$\wPIP_X$ appeared naturally, can be found in~\cite{KM1}.

We note that both semigroups $\PIP_X$ and $\wPIP_X$ admit
realisations as semigroups of \emph{difunctional binary relations}.
These special relations have been studied in a series of
works~\cite{Bre, Bre1,Sch, Sch1}. Using this
realisation $\wPIP_X$ has already appeared in~\cite{Ver}.

\section{Definitions}\label{sec:def}

Throughout the paper for a set $X$ we will denote by $'$ a bijection
from $X$ onto itself such that $(x')'=x$ for every $x\in X$.

\subsection{$\C_X$ and $\IP_X$}
First we define $\C_X$. The carrier of $\C_X$ is the set of all
partitions of $X\cup X'$ into nonempty subsets. We realise these
partitions as diagrams with two strands of vertices, top vertices
indexed by $X$ and bottom vertices indexed by $X'$. For $\alpha\in
\C_X$ two vertices of the corresponding diagram belong to the same
``connected component'' if and only if they belong to the same set
of the partition $\alpha$ (notice that there may be many different
ways of presenting an element $\alpha\in \C_X$ as a diagram, we
treat two diagrams corresponding to the same $\alpha$ as equal). The
multiplication is defined as follows: given $\alpha,\beta\in\C_X$ we
identify the bottom vertices of $\alpha$ with the corresponding top
vertices of $\beta$, which uniquely defines the connection of the
remaining vertices (which are the top vertices of $\alpha$ and the
bottom vertices of $\beta$). We set the diagram obtained in this way 
to be the product $\alpha\beta$. The formal definition of the product
$\alpha\beta$ is as follows:

Let $\alpha$, $\beta\in \C_X$, and $\equiv_{\alpha}$ and
$\equiv_{\beta}$ be the correspondent equivalence relations on
$X\cup X'$. Then the relation $\equiv_{\alpha\beta}$ is defined by:
\begin{itemize}
\item
For $i,j\in X$ we have $i\equiv_{\alpha\beta} j$ if and only if
$i\equiv_{\alpha}j$ or there exists a sequence $s_1, \dots, s_m$,
$m$ even, such that $i\equiv_{\alpha}s'_1$, $s_1\equiv_{\beta} s_2$,
$s'_2\equiv_{\alpha} s'_3$, and so on, $s_{m-1}\equiv_{\beta} s_m$,
$s'_m\equiv_{\alpha} j$.
\item
For $i,j\in X$ we have $i'\equiv_{\alpha\beta} j'$ if and only if
$i'\equiv_{\beta}j'$ or there exists a sequence $s_1, \dots, s_m$,
$m$ even, such that $i'\equiv_{\beta}s_1$, $s'_1\equiv_{\alpha}
s'_2$, $s_2\equiv_{\beta} s_3$, and so on, $s'_{m-1}\equiv_{\alpha}
s'_m$, $s_m\equiv_{\beta} j'$.
\item
For $i,j\in X$ we have $i\equiv_{\alpha\beta} j'$ if and only if
there exists a sequence $s_1, \dots, s_m$, $m$ odd, such that
$i\equiv_{\alpha}s'_1$, $s_1\equiv_{\beta} s_2$,
$s'_2\equiv_{\alpha} s'_3$, and so on,
$s'_{m-1}\equiv_{\alpha}s'_m$, $s_m\equiv_{\beta} j'$.
\end{itemize}

We will call this multiplication of partitions the {\em natural
multiplication}. An example of multiplication of elements from
$\C_8$ is given on Figure~\ref{fig:f1}.
\begin{figure}
\begin{minipage}[b]{0.45\linewidth}
$$
\xymatrix@=10pt{{\bullet}&{\bullet}\ar@{-}@/^/[r]\ar@{-}[dl]&{\bullet}\ar@{-}[d]&{\bullet}\ar@{-}[drrr]&
{\bullet}\ar@{-}@/^/[r]&{\bullet}&{\bullet}\ar@{-}@/^/[r]&{\bullet}\\
{\bullet}\ar@{-}@/_/[r]&{\bullet}\ar@{-}@/_/[r]&{\bullet}&{\bullet}\ar@{-}@/_/[r]\ar@{}[ddr]|{\times}&
{\bullet}\ar@{-}@/_/[r]&{\bullet}&{\bullet}&{\bullet}\ar@{-}[u]\ar@{-}[ul]\\
\\
{\bullet}\ar@{.}[uu]\ar@{-}[d]&{\bullet}\ar@{.}[uu]&{\bullet}\ar@{.}[uu]\ar@{-}[drr]&{\bullet}\ar@{.}[uu]\ar@{-}[d]&
{\bullet}\ar@{.}[uu]\ar@{-}@/_/[r]&{\bullet}\ar@{.}[uu]\ar@{-}@/_/[r]&{\bullet}\ar@{.}[uu]&{\bullet}\ar@{.}[uu]\ar@{-}[d]\\
{\bullet}&{\bullet}\ar@{-}@/_/[r]&{\bullet}&{\bullet}&{\bullet}&{\bullet}\ar@{-}@/_/[r]&{\bullet}&{\bullet}\\
\\
{\bullet}&{\bullet}\ar@{-}[dl]\ar@{-}@/^/[r]&{\bullet}\ar@{-}[drr]&{\bullet}\ar@{-}[d]\ar@{}[uur]|{=}&
{\bullet}\ar@{-}@/^/[r]&{\bullet}&{\bullet}\ar@{-}@/^/[r]&{\bullet}&\\
{\bullet}\ar@{-}@/^1pc/[rrrr]&{\bullet}\ar@{-}@/^/[r]&{\bullet}&{\bullet}&
{\bullet}&{\bullet}\ar@{-}@/^/[r]&{\bullet}&{\bullet}\ar@{-}[u]\ar@{-}[ul]
}
$$
\caption{Elements of $\C_8$ and their\newline
multiplication.}\label{fig:f1}
\end{minipage}
\hspace{0.7cm} % To get a little bit of space between the figures
\begin{minipage}[b]{0.45\linewidth}
$$
\xymatrix@=10pt{
{\bullet}\ar@{-}@/^/[r]\ar@{-}[d]&{\bullet}\ar@{-}[ld]&{\bullet}\ar@{-}@/^/[r]\ar@{-}[ld]&{\bullet}\ar@{-}[d]
&{\bullet}\ar@{-}[rd]&{\bullet}\ar@{-}[ld]\ar@{-}@/^/[r]&{\bullet}\ar@{-}[rd]&{\bullet}\ar@{-}[ld]\\
{\bullet}&{\bullet}\ar@{-}@/_/[r]&{\bullet}\ar@{-}@/_/[r]&{\bullet}\ar@{}[ddr]|{\times}
&{\bullet}\ar@{-}@/_1pc/[rrr]&{\bullet}&{\bullet}&{\bullet}\\
\\
{\bullet}\ar@{.}[uu]\ar@{-}@/^/[r]\ar@{-}[rd]&{\bullet}\ar@{.}[uu]\ar@{-}[d]&{\bullet}\ar@{.}[uu]\ar@{-}[rd]&{\bullet}\ar@{-}[llld]\ar@{-}[ld]\ar@{.}[uu]
&{\bullet}\ar@{-}[d]\ar@{.}[uu]\ar@{-}@/^/[r]&{\bullet}\ar@{-}[dl]\ar@{.}[uu]&{\bullet}\ar@{-}[dl]\ar@{-}[dr]\ar@{.}[uu]&{\bullet}\ar@{-}[dl]\ar@{.}[uu]\\
{\bullet}\ar@{-}@/_/[rr]&{\bullet}&{\bullet}&{\bullet}
&{\bullet}&{\bullet}\ar@{-}@/_/[rr]&{\bullet}&{\bullet}\\
\\
{\bullet}\ar@{-}[d]\ar@{-}@/_/[r]&{\bullet}\ar@{-}@/_/[r]&{\bullet}\ar@{-}@/_/[r]&{\bullet}\ar@{-}[d]\ar@{}[uur]|{=}
&{\bullet}\ar@{-}[d]\ar@{-}@/_/[r]&{\bullet}\ar@{-}@/_/[r]&{\bullet}\ar@{-}[d]&{\bullet}\ar@{-}[lld]\ar@{-}[d]\\
{\bullet}\ar@{-}@/_/[r]&{\bullet}\ar@{-}@/_/[r]&{\bullet}\ar@{-}@/_/[r]&{\bullet}
&{\bullet}\ar@{-}@/^/[rr]&{\bullet}\ar@{-}@/_/[rr]&{\bullet}&{\bullet}
}
$$
\caption{Elements of $\IP_8$ and their\newline
multiplication.}\label{fig:f2}
\end{minipage}
\end{figure}

A one-element subset of $X\cup X'$ will be called a {\em point}, and a subset $A$ intersecting with both $X$ and $X'$
--- a {\em generalised line}. A generalised line $A$ will be called a {\em line} if $|A|=2$. By $\IP_X$ we denote the subsemigroup of $\C_X$ whose elements contain only generalised lines. On Figure~\ref{fig:f2} we give an example of multiplication of the elements of $\IP_8$.

\subsection{$\PIP_X$}

Let $\PIP_X$ be the set of all partitions of the set $X\cup X'$ into
subsets being either points or generalised lines. The set $\PIP_X$
is not closed under the natural multiplication of $\C_X$ as the
example on Figure~\ref{fig:C-PIP} shows.
\begin{figure}\begin{minipage}[b]{0.45\linewidth}
$$
\xymatrix@=10pt{{\bullet}&{\bullet}\ar@{-}@/^/[r]\ar@{-}[dl]&{\bullet}\ar@{-}[dl]&{\bullet}\ar@{-}[dl]&
{\bullet}\ar@{-}@/^/[r]\ar@{-}[d]&{\bullet}\ar@{-}[dl]&{\bullet}\ar@{-}[dl]\ar@{-}@/^/[r]&{\bullet}\ar@{-}[d]\\
{\bullet}\ar@{-}@/_/[r]&{\bullet}&{\bullet}&{\bullet}\ar@{}[ddr]|{\times}&
{\bullet}&{\bullet}\ar@{-}@/_/[r]&{\bullet}\ar@{-}@/_/[r]&{\bullet}\\
\\
{\bullet}\ar@{.}[uu]&{\bullet}\ar@{.}[uu]&{\bullet}\ar@{.}[uu]\ar@{-}[d]\ar@{-}@/_/[rr]&{\bullet}\ar@{.}[uu]&
{\bullet}\ar@{.}[uu]\ar@{-}[dl]&{\bullet}\ar@{.}[uu]\ar@{-}@/^/[r]\ar@{-}[dl]&{\bullet}\ar@{.}[uu]\ar@{-}[dll]&{\bullet}\ar@{-}[dll]\ar@{.}[uu]\ar@{-}[d]\\
{\bullet}&{\bullet}&{\bullet}\ar@{-}@/_/[r]&{\bullet}&
{\bullet}&{\bullet}\ar@{-}@/^/[rr]&{\bullet}&{\bullet}\\
\\
{\bullet}&{\bullet}\ar@{-}@/^/[r]&{\bullet}&{\bullet}\ar@{-}[dl]\ar@{-}@/^/[r]\ar@{}[uur]|{=}&
{\bullet}\ar@{-}@/^/[r]&{\bullet}\ar@{-}[dll]&{\bullet}\ar@{-}[dll]\ar@{-}@/^/[r]&{\bullet}\ar@{-}[d]&\\
{\bullet}&{\bullet}&{\bullet}\ar@{-}@/_/[r]&{\bullet}&
{\bullet}\ar@{-}@/_/[r]&{\bullet}\ar@{-}@/_/[rr]&{\bullet}&{\bullet}
}
$$
\caption{$\PIP_8$ is not closed under the natural product.}\label{fig:C-PIP}
\end{minipage}
\hspace{0.7cm} % To get a little bit of space between the figures
\begin{minipage}[b]{0.45\linewidth}
$$
\xymatrix@=10pt{{\bullet}&{\bullet}\ar@{-}@/^/[r]\ar@{-}[dl]&{\bullet}\ar@{-}[dl]&{\bullet}\ar@{-}[dl]&
{\bullet}\ar@{-}@/^/[r]\ar@{-}[d]&{\bullet}\ar@{-}[dl]&{\bullet}\ar@{-}[dl]\ar@{-}@/^/[r]&{\bullet}\ar@{-}[d]\\
{\bullet}\ar@{-}@/_/[r]&{\bullet}&{\bullet}&{\bullet}\ar@{}[ddr]|{\times}&
{\bullet}&{\bullet}\ar@{-}@/_/[r]&{\bullet}\ar@{-}@/_/[r]&{\bullet}\\
\\
{\bullet}\ar@{.}[uu]\ar@{-}[d]\ar@{-}[dr]&{\bullet}\ar@{.}[uu]&{\bullet}\ar@{.}[uu]\ar@{-}[d]\ar@{-}@/_/[rr]&{\bullet}\ar@{.}[uu]&
{\bullet}\ar@{.}[uu]\ar@{-}[dl]&{\bullet}\ar@{.}[uu]\ar@{-}@/^/[r]\ar@{-}[dl]&{\bullet}\ar@{.}[uu]\ar@{-}[dll]&{\bullet}\ar@{-}[dll]\ar@{.}[uu]\ar@{-}[d]\\
{\bullet}\ar@{-}@/_/[r]&{\bullet}&{\bullet}\ar@{-}@/_/[r]&{\bullet}&
{\bullet}&{\bullet}\ar@{-}@/^/[rr]&{\bullet}&{\bullet}\\
\\
{\bullet}&{\bullet}&{\bullet}&{\bullet}\ar@{-}[dl]\ar@{-}@/^/[r]\ar@{}[uur]|{=}&
{\bullet}\ar@{-}@/^/[r]&{\bullet}\ar@{-}[dll]&{\bullet}\ar@{-}[dll]\ar@{-}@/^/[r]&{\bullet}\ar@{-}[d]&\\
{\bullet}&{\bullet}&{\bullet}\ar@{-}@/_/[r]&{\bullet}&
{\bullet}\ar@{-}@/_/[r]&{\bullet}\ar@{-}@/_/[rr]&{\bullet}&{\bullet}
}
$$
\caption{Elements of $\PIP_8$ and their
multiplication.}\label{fig:f3}
\end{minipage}

\end{figure}

However, we can define an associative multiplication on $\PIP_X$ as follows. Let $\overline{x}\notin X$. For every $\alpha\in\PIP_X$ set ${\overline{\alpha}}\in\IP_{X\cup\{\overline{x}\}}$ to be the element such that its blocks are the blocks of $\alpha$ plus one more block consisting of $\overline{x}$, $\overline{x}'$ and all points
of $\alpha$. Denote by $\varphi$ the injection, which maps $\alpha\in \PIP_X$ to ${\overline{\alpha}}\in\IP_{X\cup\{\overline{x}\}}$. Observe that $\gamma\in\IP_{X\cup\{\overline{x}\}}$ belongs to the image of $\varphi$ if and only if $\overline{x}\equiv_{\gamma}\overline{x}'$.
This enables us to define an associative multiplication on $\PIP_X$
as follows:
$$
\alpha\star\beta=\varphi^{-1}({\overline{\alpha}}{\overline{\beta}}).
$$

In terms of the diagrams we have the following interpretation of the
operation $\star$. Connect the bottom vertices of $\alpha$ with the
top vertices of $\beta$. Then two elements $a,b$ from the union of
the top vertices of $\alpha$ and the bottom ones of $\beta$ belong
to the same block of $\alpha\star\beta$ if and only if $a=b$, or $a$
and $b$ are connected and neither of them is connected to a point.
On Figure~\ref{fig:f3} we give an example of multiplication of the
elements from $\PIP_8$.

\subsection{$\wPIP_X$}
There is another way to define a multiplication on the set $\PIP_X$.
Given $\alpha,\beta$ from the set $\PIP_X$, there is a unique
element $\gamma=\alpha\circ\beta\in\PIP_X$ such that for $i,j\in X$,
$i\equiv_{\gamma}j'$ if and only if $i$ belongs to some generalised
line $A$ of $\alpha$ and $j'$ belongs to some generalised line $B$
of $\beta$ such that $A\cap X'=(B\cap X)'$. We give an example of
multiplication of elements from the set $\PIP_X$ in this way on
Figure~\ref{fig:f4}. It is easy to see that $\circ$ gives rise to a
semigroup $\wPIP_X$ on the set $\PIP_X$.

Observe, that while being closed under $\star$, $\IP_X$ is not
closed under $\circ$, which is illustrated on
Figure~\ref{fig:IP-notin-wPIP}. Besides, the $\circ$-product of the
two elements of $\PIP_8$ from Figure~\ref{fig:f3} is the element,
all the blocks of which are points. This element is a zero with
respect to both $\star$ and $\circ$. In the sequel we will denote
this element just by $0$.

\begin{figure}\begin{minipage}[b]{0.45\linewidth}
$$
\xymatrix@=10pt{{\bullet}&{\bullet}\ar@{-}@/^/[r]\ar@{-}[dl]&{\bullet}\ar@{-}[dl]&{\bullet}\ar@{-}[dl]&
{\bullet}\ar@{-}@/^/[r]\ar@{-}[d]&{\bullet}\ar@{-}[dl]&{\bullet}\ar@{-}[dl]\ar@{-}@/^/[r]&{\bullet}\ar@{-}[d]\\
{\bullet}\ar@{-}@/_/[r]&{\bullet}&{\bullet}&{\bullet}\ar@{}[ddr]|{\times}&
{\bullet}&{\bullet}\ar@{-}@/_/[r]&{\bullet}\ar@{-}@/_/[r]&{\bullet}\\
\\
{\bullet}\ar@{.}[uu]\ar@{-}[d]\ar@{-}@/^/[r]&{\bullet}\ar@{.}[uu]\ar@{-}[dl]&
{\bullet}\ar@{.}[uu]\ar@{-}[d]\ar@{-}@/_/[rr]&{\bullet}\ar@{.}[uu]&
{\bullet}\ar@{.}[uu]\ar@{-}[dl]&{\bullet}\ar@{.}[uu]\ar@{-}[dl]&{\bullet}\ar@{.}[uu]&{\bullet}\ar@{-}[dll]\ar@{.}[uu]\ar@{-}[d]\\
{\bullet}&{\bullet}&{\bullet}\ar@{-}@/_/[r]&{\bullet}&
{\bullet}&{\bullet}\ar@{-}@/^/[rr]&{\bullet}&{\bullet}\\
\\
{\bullet}&{\bullet}\ar@{-}@/^/[r]\ar@{-}[dl]&{\bullet}\ar@{-}[dll]&{\bullet}\ar@{}[uur]|{=}&
{\bullet}&{\bullet}&{\bullet}&{\bullet}\\
{\bullet}&{\bullet}&{\bullet}&{\bullet}&
{\bullet}&{\bullet}&{\bullet}&{\bullet}
}
$$
\caption{Elements of $\wPIP_8$ and their
multiplication.}\label{fig:f4}
\end{minipage}
\hspace{0.7cm} % To get a little bit of space between the figures
\begin{minipage}[b]{0.45\linewidth}
$$
\xymatrix@=10pt{
{\bullet}\ar@{-}@/^/[r]\ar@{-}[d]&{\bullet}\ar@{-}[ld]&{\bullet}\ar@{-}@/^/[r]\ar@{-}[ld]&{\bullet}\ar@{-}[d]
&{\bullet}\ar@{-}[rd]&{\bullet}\ar@{-}[ld]\ar@{-}@/^/[r]&{\bullet}\ar@{-}[dll]&{\bullet}\ar@{-}[dll]\ar@{-}[d]\\
{\bullet}&{\bullet}\ar@{-}@/_/[r]&{\bullet}\ar@{-}@/_/[r]&{\bullet}\ar@{}[ddr]|{\times}
&{\bullet}&{\bullet}\ar@{-}@/_/[rr]&{\bullet}&{\bullet}\\
\\
{\bullet}\ar@{.}[uu]\ar@{-}[d]\ar@{-}[drr]&{\bullet}\ar@{.}[uu]\ar@{-}[d]\ar@{-}@/^/[r]&
{\bullet}\ar@{.}[uu]\ar@{-}[dl]&{\bullet}\ar@{.}[uu]\ar@{-}[d]\ar@{-}@/^/[rrr]&
{\bullet}\ar@{.}[uu]\ar@{-}[d]\ar@{-}[drrr]&{\bullet}\ar@{.}[uu]\ar@{-}@/^/[rr]\ar@{-}[dr]&
{\bullet}\ar@{.}[uu]\ar@{-}[dlll]&{\bullet}\ar@{.}[uu]\ar@{-}[dl]\\
{\bullet}\ar@{-}@/_/[rr]&{\bullet}&{\bullet}&{\bullet}&
{\bullet}\ar@{-}@/_/[r]&{\bullet}\ar@{-}@/_/[rr]&{\bullet}&{\bullet}\\
\\
{\bullet}\ar@{-}@/^/[r]\ar@{-}[d]&{\bullet}\ar@{-}[dr]&{\bullet}&{\bullet}\ar@{}[uur]|{=}&
{\bullet}&{\bullet}\ar@{-}@/^/[r]\ar@{-}[dl]&{\bullet}\ar@{-}[dr]&{\bullet}\ar@{-}[dl]\\
{\bullet}\ar@{-}@/^/[rr]&{\bullet}&{\bullet}&{\bullet}&
{\bullet}\ar@{-}@/^/[r]&{\bullet}\ar@{-}@/_/[rr]&{\bullet}&{\bullet}
}
$$
\caption{$\IP_8$ is not closed under the operation
$\circ$.}\label{fig:IP-notin-wPIP}\end{minipage}
\end{figure}

%%%

In what follows we will use the following notation. Let
$\alpha\in\PIP_X$ be the element whose generalised lines are
$\{(A_i\cup B'_i)\}_{i\in I}$. Since $\alpha$ is uniquely defined by
its generalised lines, we will write $\alpha=\{(A_i\cup
B'_i)\}_{i\in I}$. We also set $\rank(\alpha)=|I|$, $\codom(\alpha) =\bigl\{ t \mid t\in
X\setminus \bigcup_{i\in I}A_i\bigr\}$, $\coran(\alpha) = \bigl\{t'
\mid t\in X\setminus \bigcup_{i\in I}B_i\bigr\}$, $\dom(\alpha)$ to be the partition  $\bigcup_{i\in I} A_i$ of the set
$X\setminus\codom(\alpha)$, $\ran(\alpha)$ --- the partition $\bigcup_{i\in I}B_i'$ of the set
$X'\setminus\coran(\alpha)$.

\section{$\PIP_X$ and
$\wPIP_X$ are inverse semigroups}\label{sec:basic_prop}

For a semigroup $S$ by $E(S)$ we denote the set of idempotents of $S$.

\begin{proposition}
$\PIP_X$ and $\wPIP_X$ are inverse semigroups.
\end{proposition}

\begin{proof}
It is sufficient to prove that the semigroups are regular and idempotents commute (see~\cite[Theorem II.1.2, p.78]{Petrich}). First we observe that idempotents in $\PIP_X$ and $\wPIP_X$ are of the
form $\bigl\{(E_i\cup E_i')\bigr\}_{i\in I}$. It follows that  both $E(\PIP_X)$ and $E(\wPIP_X)$ are semilattices.

It remains to show that $\PIP_X$ and $\wPIP_X$ are regular. Let
$\alpha=\bigl\{(A_i\cup B_i')\bigr\}_{i\in I}\in \PIP_X.$ Set
$\alpha^{-1}=\bigl\{(B_i\cup A_i')\bigr\}_{i\in I}$. Then we have
$\alpha\star\alpha^{-1}\star\alpha=\alpha$,
$\alpha^{-1}\star\alpha\star\alpha^{-1}=\alpha^{-1}$ and
$\alpha\circ\alpha^{-1}\circ\alpha=\alpha$,
$\alpha^{-1}\circ\alpha\circ\alpha^{-1}=\alpha^{-1}$.
\end{proof}

We will call the cardinality of the set of all generalised lines in
$s\in\PIP_X$ the {\em rank} of $s$ and denote it by $\rank(s)$. The
following proposition describing the structure of the Green's
relations on our semigroups is a routine to check.

\begin{proposition}\label{pr:Green_and_subgroups}
Let $a,b$ be from $\PIP_X$ or from $\wPIP_X$.
\begin{enumerate}[(1)]
\item\label{1}
$a\GR b$ if and only if $\dom(a)=\dom(b)$.
\item\label{2}
$a\GL b$ if and only if $\ran(a)=\ran(b)$.
\item\label{4}
$a\GD b$ if and only if $a\GJ b$ if and only if $\rank(a)=\rank(b)$.
\item\label{6}
All the ideals of $\PIP_X$ $($respectively $\wPIP_X)$ have the form
$$J_{\xi}=\bigl\{\alpha\in\PIP_X:~\rank(\alpha)<\xi\bigr\}$$ for certain
cardinal $\xi\leq\lmod X\rmod'$, where $\lmod X\rmod'$ is the
successor cardinal of $\lmod X\rmod$.
\end{enumerate}
\end{proposition}

\section{Fundamentality}\label{sec:fund}

Recall that an inverse semigroup $S$ is said to be \emph{fundamental} if
the \emph{maximal idempotent-separating congruence}
$$
\mu=\bigl\{(a,b)\in S\times S:~a^{-1}ea=b^{-1}eb {\text{ for all
}}e\in E(S)\bigr\}
$$
is trivial. It is well-known that $\mu$ is the largest
congruence contained in $\GH$. For $x\in X$ set
$\alpha_x=\bigl\{\{t\}\cup \{t'\}\bigr\}_{t\in X\setminus\{x\}}$.

\begin{proposition}\label{pr:fund}
Let $X$ be non-singleton. Then $\PIP_X$ and $\wPIP_X$ are
fundamental.
\end{proposition}

\begin{proof}
We will prove the statement for $\PIP_X$; for $\wPIP_X$ the proof is
similar. Suppose $(a,b)\in\mu$ for some $a,b\in\PIP_X$. Since $a\GH
b$, there are two collections of pairwise disjoint sets $A_i$, $i\in
I$, $B_i$, $i\in I$, such that
\begin{equation*}
a = \bigl\{(A_i\cup B_i')\bigr\}_{i\in I},\quad b = \bigl\{(A_i\cup
B_{\pi(i)}')\bigr\}_{i\in I}
\end{equation*}
for some bijection $\pi:~I\to I$. Let $i\in I$ and  $u_i\in A_i$.
Then  $\bigl(\alpha_{u_i}\star a,$ $\alpha_{u_i}\star b\bigr)\in\mu$
and so $(\alpha_{u_i}\star a)\GH (\alpha_{u_i}\star b)$. On the
other hand $\coran(\alpha_{u_i}\star a)=\coran(a)\cup B_i'$ and
$\coran(\alpha_{u_i}\star b)=\coran(a)\cup B_{\pi(i)}'$. Therefore
$B_i=B_{\pi(i)}$. Thus $\pi$ is the identity mapping. It follows
that $a=b$. 
\end{proof}

\begin{remark} Let $X$ be non-singleton. $\IP_X$ is not fundamental.
\end{remark}
\begin{proof} For $Y\subseteq X$
define the idempotent $\eta_Y=\{Y\cup Y', (X\setminus Y)\cup (X\setminus Y)'\}$. Let $x\in X$, $a=\eta_{x}$ and
$$
b=\bigl\{\{x\}\cup(X\setminus\{x\})', (X\setminus\{x\})\cup \{x'\}
\bigr\}\GH a.
$$
Observe that either $a^{-1}ea=\eta_{X}$ or $a^{-1}ea=a$, for every
$e\in E(\IP_X)$. In particular, $a^{-1}ea=a$ if and only if $e$
contains the block $\{x, x'\}$. Analogously, we have that either
$b^{-1}eb=\eta_{X}$ or  $b^{-1}eb=a$, for $e\in E(\IP_X)$. In
particular, $b^{-1}eb=a$ if and only if $e$ contains the block $\{x,
x'\}$. Therefore $(a,b)\in\mu$ which implies that $\mu$ is not
trivial.
\end{proof}

Note that $\IS_X$ is fundamental, \cite[p.215, ex.22]{Howie}.

\section{A generating set for $\PIP_n$}\label{sec:generating-systems}

In the case when  $X$ is $n$-set we assume that $X={\mathcal N}=\{1,2,\dots, n\}$ and in the notation for our semigroups replace lower index $X$ by $n$.

In the following sections we will need to use some generating sets
for $\PIP_n$.

Let $x,y,z\in X$ be pairwise distinct. Set
\begin{equation*}
\gamma_{x,y}=\bigl\{\{x,y\}\cup \{x'\}, \{\{t\}\cup \{t'\}\}_{t\in
X\setminus\{x,y,z\}}\bigr\};
\end{equation*}
\begin{equation*}
\xi_{x,y,z}=\bigl\{\{x,y\}\cup \{x'\}, \{z\}\cup \{y',z'\}, \{\{t\}\cup \{t'\}\}_{t\in
X\setminus\{x,y,z\}}\bigr\};
\end{equation*}
\begin{equation*}
\tau_{x,y}=\{\{x,y\}\cup\{x,y\}', \{t\}\cup\{t\}'_{t\in X\setminus\{x,y\}}\}.
\end{equation*}
Notice that $\xi_{x,y,z}\in \IP_X$. The elements $\gamma_{x,y}$ and
$\xi_{x,y,z}$ satisfy the following equalities:
\begin{equation}\label{eq:gamma}
\gamma_{x,y}\gamma_{z,y}^{-1}=
\xi_{x,y,z},\,\,\gamma_{x,y}\gamma_{x,y}^{-1}=\tau_{x,y},\,\,
\gamma_{x,y}^{-1}\gamma_{x,y}=\alpha_y;
\end{equation}
\begin{equation}\label{eq:gamma-and-S_X}
g^{-1}\gamma_{x,y}g=\gamma_{g(x),g(y)}, \text{ for any } g\in\Sym_X.
\end{equation}

\begin{lemma}\label{lm:for-5-and-6}
Let $u$ be an element of $\PIP_n$ of rank $n-1$. There are $\pi,
\tau\in\Sym_n$ such that $\pi
u\tau\in\{\tau_{1,2},\alpha_{1},\xi_{1,2,3},\gamma_{1,2},\gamma_{1,2}^{-1}\}$.
\end{lemma}

\begin{proof}
It is enough to observe that every element of rank $n-1$ coincides
with some element of the form
$\tau_{x,y}\pi,\alpha_{x}\pi,\xi_{x,y,z}\pi,\gamma_{x,y}\pi,$ or
$\gamma_{x,y}^{-1}\pi$, where $x,y,z\in X$ and $\pi\in \S_X$.
\end{proof}

It is known from~\cite[Proposition 12]{Mal} that for $n\geq 3$,
$\IP_n=\langle\Sym_n,\xi_{1,2,3}\rangle$.

\begin{lemma}\label{lm:first-g-s}
Let $n\geq 3$. Then
$\PIP_n=\langle\Sym_n,\gamma_{1,2},\gamma_{1,2}^{-1}\rangle$.
\end{lemma}

\begin{proof}
Let $a\in\PIP_n$. Consider four possible cases.

\emph{Case 1}. Suppose $a\in \IP_n$. Then from~\cite[Proposition
12]{Mal},~\eqref{eq:gamma} and \eqref{eq:gamma-and-S_X} it follows
that $a\in\langle\Sym_n,\gamma_{1,2},\gamma_{1,2}^{-1}\rangle$.

\emph{Case 2}. Suppose $a$ has a block $\{x\}$, $x\in\N$, and a
block $\{y'\}$, $y\in \N$. Let $A=\codom(a)$ and $B'=\coran(a)$.
Construct an element $q$ as follows: it contains all the generalised
lines of $a$ and, in addition, the generalised line $A\cup B'$. Then
$q\in\IP_n\subseteq\langle\Sym_n,\gamma_{1,2},\gamma_{1,2}^{-1}\rangle$.
This, $a=\alpha_{x}q\alpha_{y}$ and~\eqref{eq:gamma} imply
$a\in\langle\Sym_n,\gamma_{1,2},\gamma_{1,2}^{-1}\rangle$.

\emph{Case 3}. Suppose $a$ has a block $\{x'\}$, $x\in\N$, and has
no blocks $\{y\}$, $y\in \N$. Then there exists a generalised line
$A\cup B'$ in $a$ such that $\lmod A\rmod\geq 2$. Fix $i,j\in A$.
Set $M'=\coran(a)\ne\varnothing$. Construct the element $p$ as
follows: it contains the blocks $\{j\}\cup M'$,
$(A\setminus\{j\})\cup B'$ and all the other blocks of $p$ are all
the generalised lines of $a$ except $A\cup B'$. By the construction,
$p\in\IP_n$. Moreover, $\gamma_{i,j}p=a$. From what we have proved
in the first case now follows
$p\in\langle\Sym_n,\gamma_{1,2},\gamma_{1,2}^{-1}\rangle$, which
implies $a\in\langle\Sym_n,\gamma_{1,2},\gamma_{1,2}^{-1}\rangle$.

\emph{Case 4}. $a$ has a block $\{x\}$, $x\in\N$, and has no
blocks $\{y'\}$, $y\in \N$. This case is dual to Case 3.
\end{proof}

\begin{lemma}\label{lm:irreducible}
$\gamma_{1,2}^{-1}\notin\langle\Sym_n,\gamma_{1,2},\tau_{1,2},\xi_{1,2,3},\alpha_1\rangle$.
\end{lemma}

\begin{proof}
Assume that there are elements $a_1,\ldots,a_k$ in
$\Sym_n\{\gamma_{1,2},\tau_{1,2}, \xi_{1,2,3},\alpha_1\}\Sym_{n}$
such that $\gamma_{1,2}^{-1}=a_1\cdot\ldots\cdot a_k$. Since
$\coran(\gamma_{1,2}^{-1})=\varnothing$, it follows that
$\coran(a_k)=\varnothing$. Thus $a_k\in
\Sym_n\{\tau_{1,2},\xi_{1,2,3}\}\Sym_{n}\subseteq\IP_n$. This, in
turn, gives $\coran(a_{k-1})=\varnothing$, whereas
$a_{k-1}\in\IP_n$. Then $a_i\in\IP_n$ for all $i\leq k$ by
induction. Therefore $\gamma_{1,2}^{-1}=a_1\cdot\ldots\cdot
a_k\in\IP_n$. This is a contradiction, which completes the proof.
\end{proof}

\begin{theorem}\label{th:g-s}
Let $n\geq 3$.
\begin{enumerate}[1)]
\item $\PIP_n$ as an inverse semigroup is generated by  $\Sym_n$ and $\gamma_{1,2}$.
\item $\PIP_n$ is generated (as an inverse semigroup) by $\Sym_n$ and some $u\in\PIP_n\setminus\Sym_n$ if and only if
$u\in\Sym_n\{\gamma_{1,2},\gamma_{1,2}^{-1}\}\Sym_n$.
\end{enumerate}
\end{theorem}

\begin{proof}
The first claim follows from Lemma~\ref{lm:first-g-s}. To prove the
second one it suffices to show that
$\PIP_n=\langle\Sym_n,u,u^{-1}\rangle$ implies
$u\in\Sym_n\{\gamma_{1,2},\gamma_{1,2}^{-1}\}\Sym_n$. Let
$\PIP_n=\langle\Sym_n,u,u^{-1}\rangle$. Then $u$ is of rank $n-1$.
From  Lemma~\ref{lm:for-5-and-6} we have
$u\in\Sym_n\{\tau_{1,2},\alpha_{1},\xi_{1,2,3},\gamma_{1,2},\gamma_{1,2}^{-1}\}\Sym_n$.
Observe that $u\notin\Sym_n\{\alpha_1\}\Sym_n$ since otherwise we
would have $\langle\Sym_n,u,u^{-1}\rangle\subseteq\IS_n$, and
$u\notin\Sym_n\{\tau_{1,2},\xi_{1,2,3}\}\Sym_n$ since otherwise we
would have $\langle\Sym_n,u,u^{-1}\rangle\subseteq\IP_n$. The
statement follows.
\end{proof}

The situation with the generating sets for $\wPIP_n$ is much more complicated: one can show that 
$\wPIP_n$ can not be generated by adding to $\Sym_n$ some natural and `compact' set of elements.

\section{Maximal and maximal inverse subsemigroups}\label{sec:max}

\begin{theorem}\label{th:max-subsem}
Maximal subsemigroups of $\PIP_n$ are exhausted by the following
list:
\begin{enumerate}[1)]
\item\label{1'}
$\Sym_n\cup
J_{n-1}\cup\Sym_n\{\tau_{1,2},\alpha_1,\gamma_{1,2},\xi_{1,2,3}\}\Sym_n$;
\item\label{2'}
$\Sym_n\cup J_{n-1}\cup\Sym_n\{\tau_{1,2},\alpha_1,\gamma_{1,2}^{-1},\xi_{1,2,3}\}\Sym_n$;
\item\label{3'}
$G\cup J_n$, where $G$ runs through the set of all maximal subgroups of $\Sym_n$.
\end{enumerate}
Maximal inverse subsemigroups of $\PIP_n$ are exhausted by the
following list:
\begin{enumerate}[1)]
\item
$\Sym_n\cup J_{n-1}\cup\Sym_n\{\tau_{1,2},
\alpha_1,\xi_{1,2,3}\}\Sym_n = \langle \IP_n, \IS_n\rangle$,
\item
$G\cup J_n$, where $G$ runs through the set of all maximal
subgroups of $\Sym_n$.
\end{enumerate}
\end{theorem}

\begin{proof}
That the semigroups listed in items 1) and 2) are maximal follows
from Lemma~\ref{lm:for-5-and-6}, Lemma~\ref{lm:first-g-s} and
Lemma~\ref{lm:irreducible}. That the semigroups given in item 3) are
maximal is obvious.

Let $T$ be a maximal subsemigroup of $\PIP_n$. Then
$J_{n-1}\subseteq T$ and $G\subseteq T$, where $G$ is either $\S_n$
or a maximal subgroup of $\S_n$. If $G\neq \S_n$ then $T\subseteq
S$, where $S$ is one of the semigroups listed in item 3). Since both
$T$ and $S$ are maximal, it follows that $T=S$. Let $G=\S_n$.
Observe that we can not have $\gamma_{1,2}\in T$ and
$\gamma_{1,2}^{-1}\in T$, since otherwise we would have $T=\PIP_n$
by Lemma~\ref{lm:first-g-s}. Suppose $\gamma_{1,2}\in T$ and
$\gamma_{1,2}^{-1}\not\in T$. Then $T\subseteq S$, where
$S=J_{n-1}\cup\Sym_n\{\tau_{1,2},\alpha_1,\gamma_{1,2},\xi_{1,2,3}\}\Sym_n$.
Since both $T$ and $S$ are maximal, it follows that $T=S$. The case
$\gamma_{1,2}^{-1}\in T$ and $\gamma_{1,2}\not\in T$ is treated
similarly.

The proof of the claim about maximal inverse subsemigroups is analogous
and is left to the reader.
\end{proof}

\section{Congruences on $\IP_n$}\label{sec:congr_IP}

Let $S$ be an inverse semigroup and $E=E(S)$. We recall the definitions from~\cite[p.~118]{Petrich}. A subsemigroup $K$ of $S$ is said to be a \emph{normal subsemigroup} of $S$ if $E\subseteq
K$ and $s^{-1}Ks\subseteq K$ for all $s\in S$. A congruence
$\Lambda$ on $E$ is said to be \emph{normal} provided that for all
$e,f\in E$ and $s\in S$, $e\Lambda f$ implies $s^{-1}es\Lambda
s^{-1}fs$. The pair $(K,\Lambda)$ is said to be a \emph{congruence
pair} of $S$ if $K$ is a normal subsemigroup of $S$, $\Lambda$ is a
normal congruence on $E$ and
\begin{itemize}
\item
$ae\in K$, $e\Lambda a^{-1}a$ imply $a\in K$ for all $a\in S$ and
$e\in E$;
\item
$k\in K$ implies $kk^{-1}\Lambda k^{-1}k$.
\end{itemize}
For congruence pair $(K,\Lambda)$ of $S$ define the relation
$\rho_{(K,\Lambda)}$:
\begin{equation*}
\bigl(a\rho_{(K,\Lambda)}b\bigr)\Leftrightarrow \bigl(a^{-1}a\Lambda
b^{-1}b \text{ and } ab^{-1}\in K\bigr).
\end{equation*}
It is known (see~\cite[Theorem III.1.5, p.119]{Petrich}) that
$\rho_{(K,\Lambda)}$ is a congruence on $S$, and every congruence on
$S$ is of the form $\rho_{(K,\Lambda)}$, where $(K,\Lambda)$ is a
congruence pair of $S$.

In this section we describe all normal congruences, all normal
subsemigroups and all congruence pairs on $\IP_n$. Set $E_n=E(\IP_n)$.

\begin{lemma}\label{lm:lemma-normal-cong}
Let $e,f\in E_n$ be such that $\rank(f)\leq\rank(e)$. Then there
exists $s\in\IP_n$ such that $s^{-1}es=f$.
\end{lemma}
\begin{proof}
Suppose $e=\bigl\{E_1\cup E_1',\ldots,E_k\cup E_k'\bigr\}$,
$f=\bigl\{F_1\cup F_1',\ldots,F_l\cup F_l'\bigr\}$ where $k\geq l$.
Then $f=s^{-1}es$ for $s=\bigl\{E_1\cup F_1',\ldots,E_{l-1}\cup
F_{l-1}',(\bigcup_{i=l}^{k}E_i)\cup F_l'\bigr\}$.
\end{proof}

Let $Y\subseteq \N$. Let $\tau_Y=\{Y\cup Y', \{t\}\cup\{t\}'_{t\in \N\setminus Y}\}$. Observe that $\tau_\N=\N\cup \N'$ is the zero of $\IP_n$,
$\rank(\tau_{\N})=1$ and $\tau_{\N}$ is the only element in $\IP_n$
of rank $1$. For a set $M$ let $\iota_M$ denote the identity
relation on $M$. Set also
\begin{equation*}
I_k=\bigl\{a\in \IP_n:~\rank(a)\leq k\bigr\}~\mbox{and}
\end{equation*}
\begin{equation*}
E_n^{(k)}=\bigl\{e\in E_n:~\rank(e)\leq k\bigr\}=E_n\cap I_k.
\end{equation*}

\begin{lemma}\label{lem:zero}
Let $\Lambda$ be a normal congruence on $E_n$, $e\in E_n$ and
$\rank(e)=m$. If $e\Lambda \tau_{\N}$ then $(E_n^{(m)}\times
E_n^{(m)})\subseteq \Lambda$.
\end{lemma}
\begin{proof}
Let $f\in E_n^{(m)}$. By Lemma~\ref{lm:lemma-normal-cong} there
exists $t\in \IP_n$ such that $f=t^{-1}et$. This and the definition
of a normal congruence imply $ f=t^{-1}et\Lambda
t^{-1}\tau_{\N}t=\tau_{\N}.$
\end{proof}

The following lemma characterises normal congruences on $E(\IP_n)$:

\begin{lemma}\label{lm:Lambda-IP}
Let $\Lambda$ be a normal congruence on $E_n$. Then there is $k$
such that $\Lambda=\iota_{E_n}\cup \bigl(E_{n}^{(k)}\times
E_{n}^{(k)}\bigr)$.
\end{lemma}

\begin{proof}
Suppose $\Lambda\neq\iota_{E_n}$ (otherwise we can put $k=1$). Let
$e,f\in E_n$ be such that $e\neq f$ and $e\Lambda f$. Assume
$\rank(e)\geq\rank(f)$. Then $e\Lambda ef$ and
$\rank(e)\geq\rank(ef)$. Moreover $\rank(e)>\rank(ef)$. Indeed,
otherwise we would have $\rank(f)\geq\rank(ef)=\rank(e)\geq\rank(f)$
which would imply $\rank(ef)=\rank(e)=\rank(f)$ and then $e=ef=f$, a
contradiction. Let $\rank(e)=m\geq 2$. We will show that $(E_n^{(m)}\times E_n^{(m)})\subseteq \Lambda$. Set $B=\{1,\ldots,n-m+1\}$.
We have $\rank(\tau_B)=m$. Lemma~\ref{lm:lemma-normal-cong} implies
that there is $t\in\IP_n$ such that $t^{-1}et=\tau_{B}$. Observe
that
\begin{equation}\label{eq:1}
\rank(t^{-1}eft)<m~\mbox{and}~\tau_B\Lambda t^{-1}eft.
\end{equation}
Let $u=\tau_B t^{-1}eft\tau_B=\bigl(U_i\cup U_i'\bigr)_{i\in I}$.
Then there exists $i_0\in I$ such that $B\subseteq U_{i_0}$. We also
have $\tau_B\Lambda u$. Consider two possible cases.

\emph{Case 1}. $B=U_{i_0}$. Since $\rank(u)<\rank(e)$, it follows
that there is $j\in I\setminus\{i_0\}$ such that $C=U_j\subseteq
{\N}\setminus B=\overline{B}$ and $\lmod U_j\rmod\geq 2$. Fix
$x,y\in U_j$, $x\neq y$. It follows from $u\tau_{x,y}=u$ that
$\tau_B\Lambda u=u\tau_{x,y}\Lambda\tau_{B}\tau_{x,y}$. Let now
$p,q\in \overline{B}$, $p\neq q$. There is $g\in\Sym_n$ such that
$g(i)=i$ for all $i\in B$ and $g(x)=p$, $g(y)=q$. Then
$\tau_B=g^{-1}\tau_Bg\Lambda
g^{-1}\tau_B\tau_{x,y}g=\tau_B\tau_{p,q}$. Therefore we obtain
$$
\tau_B\Lambda\prod\limits_{p,q\in\overline{B}, p\neq q} \tau_B\tau_{p,q}=\tau_B\tau_{\overline{B}}.
$$
This implies that
\begin{equation*}
\tau_{B\cup\{x\}}=\tau_B\tau_{1,x}\Lambda\tau_B\tau_{\overline{B}}\tau_{1,x}=\tau_{\N}.
\end{equation*}
Observe that $\rank(\tau_{B\cup\{x\}})=m-1$. We have
$\bigl(E_n^{(m-1)}\times E_n^{(m-1)}\bigr)\subseteq\Lambda$ by
Lemma~\ref{lem:zero}. The latter,~\eqref{eq:1} and
Lemma~\ref{lem:zero} imply $\bigl(E_n^{(m)}\times
E_n^{(m)}\bigr)\subseteq \Lambda$, as required.

\emph{Case 2}. $B$ is a proper subset of $U_{i_0}$. Take $w\in
U_{i_0}\setminus B$. We have
$$
\tau_B\Lambda u=u\tau_{B\cup\{w\}}\Lambda
\tau_B\tau_{B\cup\{w\}}=\tau_{B\cup\{w\}}.
$$
Let $j\in\overline{B}$. There is $g\in\Sym_n$ such that $g(i)=i$ for
all $i\in B$ and $g(w)=j$. Then $\tau_B=g^{-1}\tau_Bg\Lambda
g^{-1}\tau_{B\cup\{w\}}g=\tau_{B\cup\{j\}}$ and
\begin{equation*}
\tau_B\Lambda\prod\limits_{j\in\overline{B}}\tau_{B\cup\{j\}}=\tau_{\N}.
\end{equation*}
Applying Lemma~\ref{lem:zero} we obtain $\bigl(E_n^{(m)}\times
E_n^{(m)}\bigr)\subseteq \Lambda$, as required.

We have shown that $\bigl(E_n^{(m)}\times E_n^{(m)}\bigr)\subseteq
\Lambda$ whenever $e\Lambda f$ for all idempotents $e,f$ such that
$e\neq f$ and $\rank(e)=m$. Let $k\in\mathbb{N}$, $k\leq n$, be such
that there is $e\in E_n$ of rank $k$ satisfying the condition
\begin{equation}\label{condition}
e\Lambda f \text{ for some } f\in E_n, f\neq e,
\end{equation}
while there is no $e\in E_n$ with $\rank(e)\geq k$
satisfying~\eqref{condition}. It follows that
$\Lambda=\iota_{E_n}\cup \bigl(E_{n}^{(k)}\times E_{n}^{(k)}\bigr)$.
\end{proof}

Let $e\in E_n$, $k=\rank(e)$ and $A\lhd \Sym_k$. 
It is easy to see that $H_e\cong\Sym_k$. It follows that there is unique $A_e\lhd
H_e$, such that $A_e\simeq A$. Set $N_k(A)$ to be the union of all
subgroups $A_e$, where $e$ runs through all idempotents of rank $k$.
We also set $N_{n+1}(A)=\varnothing$ whenever $A\lhd \Sym_{n+1}$.

\begin{proposition}\label{pr:cong-pairs-IP}
Let $K$ be a normal subsemigroup of $\IP_n$ and $\Lambda$ a normal
congruence on $E_n$. Then $(K,\Lambda)$ is a congruence pair of
$\IP_n$ if and only if there is $k\in {\N}$ such that
$\Lambda=\iota_{E_n}\cup \bigl(E_{n}^{(k)}\times E_{n}^{(k)}\bigr)$
and $K=E_n\cup N_{k+1}(A)\cup I_{k}$ for some  $A\lhd\Sym_{k+1}$.
\end{proposition}

\begin{proof}
The sufficiency follows from Lemma~\ref{lm:Lambda-IP} and the
observation that $E_n\cup N_{k+1}(A)\cup I_{k}$, $k\in {\N}$, is a
normal subsemigroup of $\IP_n$.

Suppose $(K,\Lambda)$ is a congruence pair of $\IP_n$.
Lemma~\ref{lm:Lambda-IP} implies that there is $k\in {\N}$ such that
$\Lambda=\iota_{E_n}\cup \bigl(E_{n}^{(k)}\times E_{n}^{(k)}\bigr)$.

Assume that $k=n$. Then $\Lambda=E_n\times E_n$. In this case we
have $K=\IP_n$. Indeed, let $a\in \IP_n$. Since $E_n\subset K$, it
follows that, in particular, $a\tau_{\N}=\tau_{\N}\in K$. We also
have $\tau_{\N}\Lambda aa^{-1}$. The first condition of the
definition of a congruence pair yields $a\in K$. Hence
$\IP_n\subseteq K$.

Assume now that $k\leq n-1$. Since $a \tau_{\N}=\tau_{\N}\in
E_n\subseteq K$ and $\tau_{\N}\Lambda a^{-1}a$ for all $a\in I_k$,
it follows that $I_k\subseteq K$. Since $d\in K$ implies
$dd^{-1}\Lambda d^{-1}d$ for all $d\in K$, it follows that all
elements $d\in K$ such that $\rank(d)\geq k+1$, belong to certain
subgroups of $\IP_n$. Let $b\in K$ and $\rank(b)=m\geq k+2$. Observe
that $m\geq 3$. Show that $b$ must be an idempotent. Since $b$ is a
group element, there exists a partition ${\N}=\bigcup_{i\in I}B_i$
such that $b=\bigl\{(B_i\cup B_{\pi(i)}')_{i\in I}\bigr\}$, $\lmod
I\rmod\geq 3$, for some bijection $\pi:~I\to I$. Show that $\pi$ is
the identity transformation of $I$. Consider $\pi$ as a permutation
from $\Sym_I$. Suppose $\pi$ is not the identity map. Consider a
cycle $(i_1,i_2,\ldots,i_l)$ of $\pi$, where $i_1,\ldots,i_l\in I$
and $l\geq 2$. If $l\geq 3$ then $b\tau_{B_{i_1}\cup B_{i_2}}$ is of
rank $m-1$, is not a group element and belongs to $K$, which is
impossible. If $l=2$ consider $j\in I\setminus\{i_1,i_2\}$ and
$b\tau_{B_{i_1}\cup B_j}$. This element is again of rank $m-1$, is
not a group element and belongs to $K$, which is also impossible.
Thus $\pi$ is the identity transformation of $I$. Therefore $b$ is
an idempotent. It follows that $K\cap(\IP_n\setminus
I_{k+1})=E_n\setminus I_{k+1}$.

Fix $e,f\in E_n$ such that $\rank(e)=\rank(f)=k+1$. Set $A_e'=K\cap
H_{e}$, $A_f'=K\cap H_{f}$. Since $K$ is self-conjugate it follows
that $A_{e}'\lhd H_e$ and $A_{f}'\lhd H_f$. Take any $s\in\IP_n$
such that $s^{-1}es=f$ and $sfs^{-1}=e$ (it is easily seen that such
an element exists). Further, from $s^{-1}Ks\subseteq K$ and
$sKs^{-1}\subseteq K$, it follows that the maps $x\mapsto s^{-1}xs$
from $A_e'$ onto $A_f'$ and $y\mapsto sys^{-1}$ from $A_f'$ onto
$A_e'$ are mutually inverse bijections, whence $\lmod
A_e'\rmod=\lmod A_f'\rmod$. It follows that an element of $K$ has
rank $k+1$ if and only if it lies in $N_{k+1}(A)$ for some $A \lhd
\Sym_{k+1}$. Thus $K=E_n\cup  N_{k+1}(A)\cup I_{k}$, and the proof
is complete.
\end{proof}
For $1\leq k\leq n$ denote by $D_k$ the set of elements of $\IP_n$
of rank $k$. Let $A\lhd \Sym_{k+1}$, $1\leq k\leq n$. Let $F_k(A)$
be the relation on $D_{k+1}$ that is defined by $(x,y)\in F_k(A)$ if
and only if $x\GH y$ and $xy^{-1}\in N_{k+1}(A)$. Set
$\rho_{k,A}=\iota_{\IP_n}\cup F_k(A)\cup\bigl(I_k\times I_k\bigr)$.
The construction implies that $\rho_{k,A}$ coincides with
$\rho_{(K,\Lambda)}$, corresponding to the congruence pair
$(K,\Lambda)$, where $K=E_n\cup N_{k+1}(A)\cup I_{k}$ and
$\Lambda=\iota_{E_n}\cup \bigl(E_{n}^{(k)}\times E_{n}^{(k)}\bigr)$.

\begin{theorem}\label{th:congruences}
Let $\rho$ be a relation on $\IP_n$. Then $\rho$ is a congruence on
$\IP_n$ if and only if $\rho=\rho_{k,A}$ for some $k, 1\leq k\leq n$
and normal subgroup $A\lhd\Sym_{k+1}$.
\end{theorem}

\begin{proof}
The claim follows from Proposition~\ref{pr:cong-pairs-IP} and
from~\cite[Theorem III.1.5]{Petrich}.
\end{proof}

We note that the formulation of Theorem~\ref{th:congruences} resembles the one of the corresponding classic Liber's result~\cite{Liber} for $\IS_n$. 

\section{Congruences on $\PIP_n$ and $\wPIP_n$}\label{sec:cognr_pip}

\subsection{Congruences on $\PIP_n$}

Let $Y\subseteq X$. Set $\alpha_{Y}=\bigl\{\{t,t'\}_{t\in {X}\setminus Y}\bigr\}.$ Notice
that $\alpha_{Y}$ is an idempotent for any $Y\subseteq X$ and that
the element $0=\alpha_{X}$ is the zero element of both $\PIP_X$ and
$\wPIP_X$.

Let $\widetilde{E}_n=E(\PIP_n)$ and
$\widetilde{E}_n^{(k)}=\bigl\{e\in \widetilde{E}_n:~\rank(e)\leq
k\bigr\}=\widetilde{E}_n\cap J_{k+1}$.

\begin{lemma}\label{lm:lemma-normal-cong-PIP}
Let $e,f\in \widetilde{E}_n$ and $\rank(f)\leq\rank(e)$. Then there
exists $s\in\PIP_n$ such that $s^{-1}\star e\star s=f$.
\end{lemma}

\begin{proof}
The proof is analogous to that of Lemma~\ref{lm:lemma-normal-cong}.
\end{proof}

As an immediate consequence we obtain the following lemma.

\begin{lemma}\label{cor:trivial}
Let $\Lambda$ be a normal congruence on $(\widetilde{E}_n,\star)$.
Then $a\Lambda 0$ implies $b\Lambda 0$ for all idempotents $b\in
J_{\rank(a)+1}$.
\end{lemma}

\begin{lemma}\label{lm:IS-for-PIP}
Let $a$ and $b$ of $\IS_n$ be two idempotents with
$\rank(a)>\rank(b)$ and $\Lambda$ --- a normal congruence on
$E(\IS_n)$. Then $a$ is $\Lambda$-related to $0$.
\end{lemma}

\begin{proof}
The proof is similar to that of Lemma~\ref{lm:Lambda-IP}.
\end{proof}

\begin{lemma}\label{lm:Lambda-PIP}
Let $\Lambda$ be a normal congruence on $(\widetilde{E}_n,\star)$.
Then there is $k\in\mathcal{N}$ such that
$\Lambda=\iota_{\widetilde{E}_n}\cup
\bigl(\widetilde{E}_{n}^{(k)}\times \widetilde{E}_{n}^{(k)}\bigr)$.
\end{lemma}

\begin{proof}
Suppose $\Lambda\neq \iota_{\widetilde{E}_n}$. Take distinct
$e,f\in\widetilde{E}_n$ such that $e\neq f$ and
$m=\rank(e)\geq\rank(f)$. Show that
$\bigl(\widetilde{E}_{n}^{(m)}\times
\widetilde{E}_{n}^{(m)}\bigr)\subseteq\Lambda$. Similarly to as it
was done in the proof of Lemma~\ref{lm:Lambda-IP} we show that
$\tau_B\Lambda u$, where $B=\{1,\ldots,n-m+1\}$ and
$u\in\widetilde{E}_n$ are such that $\rank(u)<m$ and $u=\tau_B
u=u\tau_B$. Show that there exists an element of rank $m$ which is
$\Lambda$-related to $0$. Set
\begin{equation*}
d=\bigl\{B\cup\{1'\},\{k,k'\}_{k\in\N\setminus B}\bigr\}.
\end{equation*}
Consider three possible cases.

\emph{Case 1}. Suppose $u$ contains a block $C\cup C'$, where $C$
strictly contains $B$. Then $\tau_B\Lambda
u=u\tau_C\Lambda\tau_B\tau_C=\tau_C$. This and
Lemma~\ref{lm:Lambda-IP} imply $\tau_B\Lambda\tau_{\N}$. It follows
that $\alpha_{B}=\tau_B \alpha_{1}\Lambda \tau_{\N}\alpha_{1}=0$.
Since $\rank(u)\leq\rank(\alpha_{B})$ it follows from
Lemma~\ref{cor:trivial} that $u\Lambda 0$, whence $\tau_B\Lambda 0$.

\emph{Case 2}. Suppose $u$ contains a block $\{t\}$ for some $t\in
B$. Then
\begin{equation*}
\tau_B\Lambda u=\alpha_{t}u\Lambda
\alpha_{t}\tau_B=\alpha_{B}=\alpha_{1}\ldots \alpha_{n-m+1}.
\end{equation*}
Therefore
\begin{equation*}
\alpha_{B\setminus\{1\}}=\alpha_{2}\ldots
\alpha_{n-m+1}=d^{-1}\tau_B d\Lambda d^{-1}\alpha_{1}\ldots
\alpha_{n-m+1}d=\\\alpha_{1}\ldots \alpha_{n-m+1}=\alpha_{B}.
\end{equation*}
Both $\alpha_{B\setminus\{1\}}$ and $\alpha_{B}$ belong to $\IS_n$.
In addition, $\rank(\alpha_{B\setminus\{1\}})=m$ and
$\rank(\alpha_{B})=m-1$. Applying Lemma~\ref{lm:IS-for-PIP} we
obtain $\alpha_{2}\ldots\alpha_{n-m+1}\Lambda 0$.

\emph{Case 3}. Suppose $u$ contains a block $B\cup B'$. If
$u\in\IP_n$ then Lemma~\ref{lm:Lambda-IP} ensures that
$\tau_B\Lambda\tau_{\N}$. Applying the same arguments as in the
first case, we conclude that $\tau_{B}\Lambda 0$. Otherwise there is
$j\in\N\setminus B$ such that $\tau_B\Lambda \tau_B \alpha_{j}$.
Then
\begin{equation*}
\alpha_{B\setminus\{1\}}=\alpha_{2}\ldots
\alpha_{n-m+1}=d^{-1}\tau_B d\Lambda d^{-1}\tau_B
\alpha_{j}d=\alpha_{B}\alpha_{j}.
\end{equation*}
Observe that
$\rank(\alpha_{B}\alpha_{j})<\rank(\alpha_{B\setminus\{1\}})=m$.
This and Lemma~\ref{lm:IS-for-PIP} imply
$\alpha_{B\setminus\{1\}}\Lambda 0$.

Lemma~\ref{cor:trivial} implies that
$\widetilde{E}_n^{(m)}=\widetilde{E}_n\cap J_{m+1}$ lies in some
$\Lambda$-class. Applying the same arguments as at the end of the
proof of Lemma~\ref{lm:Lambda-IP}, we obtain that there is $k\in\N$
such that $\Lambda=\iota_{\widetilde{E}_n}\cup
\bigl(\widetilde{E}_{n}^{(k)}\times \widetilde{E}_{n}^{(k)}\bigr)$.
\end{proof}

For $A\lhd\Sym_k$ we construct the set $\widetilde{N}_k(A)$ and the
relation $\widetilde{F}_k(A)$ similarly to as we constructed
$N_k(A)$ and $F_k(A)$ in Section~\ref{sec:congr_IP}. Set
$\widetilde{\rho}_{k,A}=\iota_{\PIP_n}\cup
\widetilde{F}_k(A)\cup\bigl(J_{k+1}\times J_{k+1}\bigr)$. The proof
of the following statement is analogous to that of
Proposition~\ref{pr:cong-pairs-IP}.

\begin{proposition}\label{pr:cong-pairs-PIP}
Let $K$ be a normal subsemigroup of $\PIP_n$ and $\Lambda$ be a
normal congruence on $\widetilde{E}_n$. Then $(K,\Lambda)$ is a
congruence pair of $\PIP_n$ if and only if there is
$k\in\mathcal{N}$ such that  $\Lambda=\iota_{\widetilde{E}_n}\cup
\bigl(\widetilde{E}_{n}^{(k)}\times \widetilde{E}_{n}^{(k)}\bigr)$
and $K=\widetilde{E}_n\cup \widetilde{N}_{k+1}(A)\cup J_{k+1}$ for
some $A\lhd\Sym_{k+1}$.
\end{proposition}

The description of congruences on $\PIP_n$ can be formulated now in
the same way as Theorem~\ref{th:congruences}.

For the semigroup $\wPIP_n$ the arguments are similar. In
particular, we observe that an analogue of Lemma~\ref{lm:Lambda-PIP}
holds. After this, it is easy to conclude that sets of congruences
on $\PIP_n$ and $\wPIP_n$ coincide.

\section{Completely isolated subsemigroups of $\IP_n$, $\PIP_n$ and
$\wPIP_n$}\label{sec:compl_isol}

From now on suppose that $n\geq 2$. Recall that a subsemigroup $T$ of a
semigroup $S$ is called \emph{completely isolated} provided that
$ab\in T$ implies either $a\in T$ or $b\in T$ for all $a,b\in S$. A
subsemigroup $T$ of a semigroup $S$ is called {\em isolated}
provided that $a^k\in T$, $k\geq 1$, implies $a\in T$ for all $a\in
T$. A completely isolated subsemigroup is isolated, but the converse
is not true in general.

We begin this section with several general observations, which will
be needed for the sequel and are also interesting on their own.

\begin{lemma}\label{pr:very-nice}
Let $S$ be a semigroup with an identity element $1$ and the group of
units $G$. Suppose $S\setminus G$ is a subsemigroup of $S$. Then $G$
is completely isolated and the map $T\mapsto T\cup G$ is a
bijection from the set of all completely isolated subsemigroups,
which are disjoint with $G$, to the set of all completely isolated
subsemigroups, which contain $G$ as a proper subsemigroup.
\end{lemma}

\begin{proof}
Obviously, $G$ is a completely isolated subsemigroup.
Suppose that $T$ is a completely isolated subsemigroup such that
$T\cap G=\varnothing$. Observe that $T\cup G$ is a subsemigroup of
$S$. Indeed, let $g\in G$ and $t\in T$. Since $T$ is completely
isolated and disjoint with $G$, the inclusion $g^{-1}\cdot gt=t\in
T$ implies $gt\in T\subset T\cup G$. Similarly, $tg\cdot g^{-1}=t\in
T$ implies $tg\in T\subset T\cup G$. Let now $ab\in T\cup G$.
Consider two possible cases.

{\em Case 1.} Suppose $ab\in G$. Since $S\setminus G$ is a
subsemigroup of $S$, it follows that either $a\in G$ or $b\in G$.

{\em Case 2.} Suppose $ab\in T$. Since $T$ is completely isolated,
it follows that either $a\in T$ or $b\in T$.

Therefore, either $a\in G\cup T$ or
$b\in G\cup T$. Hence $T\cup G$ is completely isolated.

Now suppose that $T$ is a completely isolated subsemigroup with
$T\supset G$ and prove that $T\setminus G$ is  completely isolated
as well. Let $a,b\in T\setminus G=T\cap (S\setminus G)$. Then $ab\in
T\setminus G$ as both $T$ and $S\setminus G$ are subsemigroups of
$S$, proving that $T\setminus G$ is a semigroup. Suppose
$ab\in T\setminus G$ and show that at least one of the elements $a$,
$b$ lies in $T\setminus G$. Since $T\setminus G\subset T$ and $T$ is
completely isolated, it follows that at least one of the elements
$a,b$ belongs to $T$. Suppose $a\in T$ (the case when $b\in T$ is
treated similarly). If $a\in T\setminus G$, we are done. If $a\in G$
we have $b=a^{-1}\cdot ab\in T$. Moreover, $b\in T\setminus G$ as
the inclusion $b\in G$ would imply $ab\in G$. Hence $T\setminus G$
is completely isolated.
\end{proof}

\begin{lemma}\label{pr:isol}
Let $S$ be a semigroup, $e\in E(S)$ and $G=G(e)$ ---
the maximal subgroup of $S$ with the identity element $e$. Suppose
$G$ is periodic and $T$ is an isolated subsemigroup of $S$ such that
$T\cap G\neq \varnothing$. Then $T\supseteq G$.
\end{lemma}

\begin{proof}
Let $a\in T\cap G$. There is $m\in {\mathbb N}$ such that $a^{m}=e$,
which implies $e\in T$. Let $b\in G$. Since $G$ is periodic, $b^k=e$
for certain $k \in {\mathbb N}$. The statement follows.
\end{proof}

\begin{corollary}\label{cor:lovely}
Let $S$ be a semigroup with the group of units $G$. Suppose that
$S\setminus G$ is a subsemigroup of $S$ and that $G$ is periodic.
\begin{enumerate}
\item
If $T_i, i\in I,$ is the full list of completely isolated
subsemigroups of $S$, which are disjoint with $G$, then $T_i, i\in
I, T_i\cup G, i\in I, G$ is the full list of completely isolated
subsemigroups of $S$.
\item
If $T_i, i\in I,$ is the full list of completely isolated
subsemigroups of $S$, which contain $G$ as a proper subsemigroup,
then $T_i, i\in I, T_i\setminus G, i\in I, G$ is the full list of
completely isolated subsemigroups of $S$.
\end{enumerate}
\end{corollary}

\begin{proof}
The proof follows from Lemma~\ref{pr:very-nice} and
Lemma~\ref{pr:isol}.
\end{proof}

\subsection{Completely isolated subsemigroups of $\IP_n$}

\begin{theorem}\label{th:Completely-isolated-IP}
Let $n\geq 2$. The semigroups $\IP_n$, $\Sym_n$ and $\IP_n\setminus
\Sym_n$ and only them are completely isolated subsemigroups of the
semigroup $\IP_n$.
\end{theorem}

\begin{proof}
For $n=2$ the proof is easy. Suppose $n\geq 3$. That all the
subsemigroups given in the formulation are completely isolated
follows from the definition.

Let $T$ be a completely isolated subsemigroup of $\IP_n$ containing
$\S_n$ as a proper subsemigroup. Applying
Corollary~\ref{cor:lovely}, it is enough to prove that $T=\IP_n$.
Show that $T$ contains some element from $\Sym_n\xi_{1,2,3}\Sym_n$.
Indeed, consider $g\in T\setminus \S_n$. Due to
$\IP_n=\langle\Sym_n,\xi_{1,2,3}\rangle$ (~\cite[Proposition
12]{Mal}) we can write
\begin{equation*}
g=g_1\xi_{1,2,3}g_2\xi_{1,2,3}\cdots \xi_{1,2,3}g_{k+1},
\end{equation*}
where $k\geq 1$ and $g_1\dots, g_{k+1}\in \S_n$. If $k>1$ we have
that either $g_1\xi_{1,2,3}g_2\xi_{1,2,3}$ $\cdots \xi_{1,2,3}g_k\in
T$ or $\xi_{1,2,3}g_{k+1}\in T$, since $T$ is completely isolated.
The claim follows by induction.

Now we can assert that $\xi_{1,2,3}\in T$ as $T\supset \S_n$ by the
assumption. This together with
$\IP_n=\langle\Sym_n,\xi_{1,2,3}\rangle$ implies $T=\IP_n$.
\end{proof}

\subsection{Completely isolated subsemigroups of $\PIP_n$}

\begin{theorem}\label{th:PIP-comp-isol}
Let $n\geq 2$.  The semigroups $\PIP_n$, $\Sym_n$ and
$\PIP_n\setminus\Sym_n$ and only them are the completely isolated
subsemigroups of the semigroup $\PIP_n$.
\end{theorem}

For the
proof of Theorem~\ref{th:PIP-comp-isol} we will need two auxiliary
lemmas:

\begin{lemma}\label{lm:PIP-compl-isol}
Let $a\in\PIP_n\setminus \IP_n$. Then there are $k\geq 1$ and
$g_1,\ldots,g_k$ of $\Sym_n$ such that $ag_1ag_2a\ldots ag_ka=0$.
\end{lemma}

\begin{proof}
The statement follows from the observation that $a$ has at least one
point.
\end{proof}

\begin{lemma}\label{cor:a}
Let $T$ be a completely isolated subsemigroup of $\PIP_n$ such that
$(\{0\}\cup \S_n)\subset T$. Then $\PIP_n\setminus\IP_n\subset T$.
\end{lemma}
\begin{proof}
Let $a\in\PIP_n\setminus\IP_n$. By Lemma~\ref{lm:PIP-compl-isol} we
have $a g_1a\ldots g_ka=0$ for some $g_1,\ldots,g_k \in\Sym_n$.
Since $T$ is completely isolated, it follows that either
$ag_1ag_2\ldots ag_k\in T$ or $a\in T$. If $a\in T$ then we are
done. Otherwise, we have $ag_1ag_2\ldots ag_{k-1}a\in T$. The
statement follows by induction.
\end{proof}

\begin{proof}[Proof of Theorem~\ref{th:PIP-comp-isol}.]
That all the listed semigroups are completely isolated is checked
directly. Let $T$ be a completely isolated subsemigroup of $\PIP_n$
strictly containing $\S_n$. In view of Corollary~\ref{cor:lovely} it
is enough to show that $T=\PIP_n$.

First assume that $T\setminus \IP_n\neq \varnothing$. Take any $a\in
T\setminus\IP_n$. Since $a\in\PIP_n\setminus \IP_n$, it follows from
Lemma~\ref{lm:PIP-compl-isol} that $0\in T$. Applying
Lemma~\ref{cor:a} we obtain the inclusion $\PIP_n\setminus
\IP_n\subseteq T$.

Consider the element
\begin{equation*}
w=\gamma_{1,2}=\bigl\{\{1,2,1'\},\{t,t'\}_{t\in\N\setminus\{1,2\}}\bigr\}\in\PIP_n.
\end{equation*}
Since $w^2=(w^{-1})^2=\alpha_{1}\alpha_{2}\in\IS_n\subseteq T$, we
conclude that $w\in T$ and $w^{-1}\in T$, which implies $ww^{-1}\in
T$. From the other hand, $ww^{-1}=\tau_{{1,2}}\in\IP_n\setminus
\S_n$. It follows that $ww^{-1}\in
T\cap\bigl(\IP_n\setminus\Sym_n\bigr)$. Observe that
$\tau_{\N}\in\langle\Sym_n,\tau_{{1,2}}\rangle\subseteq T$. It is
easy to see that $T\cap\IP_n$ is a completely isolated subsemigroup
of $\IP_n$. In addition, $T\cap\IP_n$ contains $\S_n$ as a proper
subsemigroup. Applying Theorem~\ref{th:Completely-isolated-IP} we
obtain $\IP_n\subseteq T$. It follows that $T=\PIP_n$.

Assume now that $T\setminus \IP_n =\varnothing$, that is,
$T\subseteq \IP_n$. Let $a=\bigl\{(A_i\cup B_i')_{i\in I}\bigr\}\in
T\setminus \S_n$. Since $a\not\in \S_n$, there exists $j\in I$ such
that $|B_j|\geq 2$. Fix some $x\in B_j$ and consider the elements
\begin{equation*}
b=\bigl\{(A_i\cup B_i')_{i\in I\setminus\{j\}}, A_j\cup\{x'\}\bigr\} \text{ and }
\end{equation*}
\begin{equation*}
c=\bigl\{(B_i\cup B_i')_{i\in I\setminus\{j\}}, \{x\}\cup B_j'
\bigr\}
\end{equation*}
of  $\PIP_n\setminus \IP_n$. We have $bc=a\in T$ by the
construction. Therefore, $b\in T\subset \IP_n$ or $c\in T\subset
\IP_n$. We obtained a contradiction, which shows that the inclusion
$T\subseteq\IP_n$ is impossible. The proof is complete.
\end{proof}

\subsection{Completely isolated subsemigroups of $\wPIP_n$}

\begin{theorem}\label{th:wPIP-comp-isol}
Let $n\geq 2$. All completely isolated subsemigroups of $\wPIP_n$
are exhausted by the following list: $\wPIP_n$, $\Sym_n$ and
$\wPIP_n\setminus\Sym_n$.
\end{theorem}

\begin{lemma}\label{lem:aux}
Let $e\in E(\wPIP_n)\setminus \S_n$. Then there exists $a\in
\wPIP_n$ such that $aa^{-1}=e$ and $a^2=(a^{-1})^2\in \IS_n\setminus
\S_n$.
\end{lemma}

\begin{proof}
If $e\in \IS_n$ we can set $a=e$. Otherwise, let
$e=\bigl\{(A_i\cup A_i')_{i\in I}, \{t,t'\}_{t\in J} \bigr\},$
where $\N\setminus \bigl((\bigcup_{i\in I}A_i)\cup J\bigr)$ is non-empty and
$|A_i|\geq 2$, $i\in I$. Since $e\not\in \IS_n$, it follows  that
$I\neq\varnothing$. Take $x_i\in A_i$, $i\in I$. Set $a=\bigl\{(A_i\cup x_i')_{i\in I}, \{t,t'\}_{t\in J} \bigr\}.$
We have that $aa^{-1}=e$ and
$a^{2}=(a^{-1})^2=\bigl\{\{t,t'\}_{t\in J}\bigr\}\in\IS_n\setminus \S_n$.
\end{proof}

The following statement follows from Lemma~\ref{lem:aux}.

\begin{corollary}\label{cor:new}
Let $T$ be an isolated subsemigroup of $\wPIP_n$. If $\IS_n\setminus
\S_n\subseteq T$, then $\wPIP_n\setminus \S_n\subseteq T$.
\end{corollary}

We will need the following fact, see \cite[Chapter 5]{GM1}.
\begin{lemma}\label{Mazorchuk-Gan-PIP-compl-isol}
All completely isolated subsemigroups of $\IS_n$ are
exhausted by the following list: $\IS_n$, $\Sym_n$ and
$\IS_n\setminus\Sym_n$.
\end{lemma}

\begin{proof}[Proof of Theorem \ref{th:wPIP-comp-isol}]
It is straightforward to verify that $\wPIP_n$, $\Sym_n$ and
$\wPIP_n\setminus\Sym_n$ are completely isolated.
Let now $T$ be a completely isolated subsemigroup of $\wPIP_n$. If
$T\cap\Sym_n\ne\varnothing$ then $T\supset \S_n$ by
Lemma~\ref{pr:isol}. Assume that $T\setminus\Sym_n\ne\varnothing$.
It is enough to prove that $\wPIP_n\setminus\Sym_n\subseteq T$.

Let $b\in T\setminus \S_n$. There is $k$ such that $b^k=e$ is an
idempotent. Let $a\in \wPIP_n$ be such that $aa^{-1}=e$ and
$f=a^2=(a^{-1})^2\in \IS_n\setminus \S_n$ (such an element exists by
Lemma~\ref{lem:aux}). Then $f\in T\cap \IS_n$. Applying
Lemma~\ref{Mazorchuk-Gan-PIP-compl-isol} we have
$\IS_n\setminus\Sym_n\subseteq T$. Finally,
$\wPIP_n\setminus\Sym_n\subseteq T$ by Corollary~\ref{cor:new}.
\end{proof}

\section{Isolated subsemigroups of $\IP_n$, $\PIP_n$ and
$\wPIP_n$}\label{sec:isol}

\subsection{Isolated subsemigroups of $\IP_n$}

\begin{proposition}\label{pr:is_ip}
Let $e\in \IP_n$ be an idempotent of rank $n-1$, that is, $e=\tau_A$
for some $A\subset \N$ with $|A|=2$. Then $G(e)$ is an isolated
subsemigroup of $\IP_n$.
\end{proposition}
\begin{proof}
Assume that $a\in\IP_n$ is such that $a^k\in G(e)$ for some $k\geq
1$. Since $G(e)$ is finite, we can assume that $a^k=e$. We are to
show that $a\in G(e)$. Since $\rank(a^k)=n-1$, it follows that
$\rank(a)\geq n-1$. Hence $\rank(a)=n-1$. But $\rank(a)=\rank(a^k)$
implies that $a\GD a^2$, which implies that $a\GH a^2$ (since $\IP_n$ is finite), which means that $a\in
G(e)$.
\end{proof}

\begin{theorem}\label{th:is_ip}
The semigroups $\IP_n$, $\S_n$, $\IP_n\setminus \S_n$ and $G(e)$,
where $e$ is an idempotent of rank $n-1$ and only them are isolated
subsemigroups of $\IP_n$.
\end{theorem}
\begin{proof}
That all the listed subsemigroups are isolated follows from
Proposition~\ref{pr:is_ip} and
Theorem~\ref{th:Completely-isolated-IP}.

Assume that $T\neq \S_n$ is an isolated subsemigroup of $\IP_n$. Then $T\setminus \S_{n}\neq\varnothing$. Let
$a\in T\setminus \S_{n}$. Going, if necessary, to some power of $a$,
we may assume that $a$ is an idempotent. Let us show
that $T$ contains some idempotent of rank $n-1$.

Suppose first that $a$ has some block $A\cup A'$ with $A\subseteq
\N$, $|A|\geq 3$. Let $A=\{t_1,\dots t_k\}$. Consider $b\in \IP_n$
such that it contains all the blocks of $a$, except $A\cup A'$, and
instead of $A\cup A'$ it has two blocks: $\{t_1,\dots,
t_{k-1},t_1'\}$ and $\{t_k, t_2',\dots, t_k'\}$. The construction
implies $b^2=(b^{-1})^2=a$, whence $b,b^{-1}\in T$. It follows that
$bb^{-1}\in T$. This element is an idempotent, contains all the
blocks of $a$, except $A\cup A'$, and instead of $A\cup A'$ it
contains two blocks: $(A\setminus \{t_k\})\cup(A\setminus \{t_k\})'$
and $\{t_k, t_k'\}$. Applying the described procedure as many times
as needed we obtain that there $T$ contains an idempotent $e$ such
that $|A|\leq 2$ for each block $A\cup A'$, $A\subseteq \N$, of $e$.

Suppose now that $e\in E(T)$ contains two
blocks $\{t_1,t_2\}\cup \{t_1,t_2\}'$ and $\{t_3,t_4\}\cup
\{t_3,t_4\}'$, $t_1, t_2, t_3, t_4\in \N$. Let $a\in \IP_n$ be the
element whose blocks are all the blocks of $e$, except
$\{t_1,t_2\}\cup \{t_1,t_2\}'$ and $\{t_3,t_4\}\cup \{t_3,t_4\}'$,
and instead of these two blocks it contains the following three
blocks: $\{t_1,t_3'\}$, $\{t_2,t_4'\}$, $\{t_3,t_4,t_1',t_2'\}$. The
construction of $a$ implies that $a^2=(a^{-1})^2=e$, which implies
$a,a^{-1}\in T$. It follows that $aa^{-1}\in T$. Observe that $aa^{-1}\in E(T)$. 
This element contains all the blocks
of $e$, except $\{t_1,t_2\}\cup \{t_1,t_2\}'$. In addition, it has two
blocks $\{t_1,t_1'\}$ and $\{t_2,t_2'\}$. Therefore, $aa^{-1}$ has
fewer blocks of the form $A\cup A'$ with $A\subset \N$, $|A|=2$ than
$e$. Applying this procedure as many times as required we obtain
that $T$ contains some idempotent $e=\tau_A$ with $|A|=2$. Therefore, $T$ contains some idempotent $e$ of rank $n-1$.

If $e$ is
the only idempotent of $T$ we have $T=G(e)$. Suppose now that,
except $e$, $T$ has some other idempotent, say, $f$. We will show
that $\tau_{\N}\in T$. If $n=2$ this is obvious. Suppose $n\geq 3$.
In view of Lemma~\ref{pr:isol} $G(e), G(f)\subset T$. Let
$e=\tau_A$, where $A=\{t_1, t_2\}$. Consider two possible cases.

{\em Case 1.} Suppose $\rank(f)\leq n-1$. Since $f\neq e$ it follows
that $f$ has a block $B\cup B'$ with $B\subseteq \N$, $|B|\geq 2$
and $B\setminus A\neq \varnothing$. Fix some $t_3\in B\setminus A$
and $s\in B$, $s\neq t_3$. For each $i\in \N\setminus\{t_1,t_2\}$
consider the transposition $\pi_i$ of $G(e)$ which swaps $i$ and
$t_3$. Then the idempotent $e_i=(\pi_i f)(\pi_i f)^{-1}$ has a block
$C\cup C'$, $C\subseteq \N$ with $i, s\in C$. Now consider the
transposition $\pi_1\in G(e)$ which switches the blocks
$\{t_1,t_2\}$ and $\{t_3\}$. Then the idempotent $e_1=(\pi_1
f)(\pi_1 f)^{-1}$ has a block $C\cup C'$, $C\subseteq \N$, with
$t_1,t_2, s\in C$. The product of all the constructed idempotents
$e_i$, $i\in \N\setminus\{t_2\}$, equals $\tau_{\N}$.

{\em Case 2.} Suppose $\rank(f)=n$, that is, $f=1$.
Then $\S_n\subseteq T$. Conjugating $e$ by each transposition of
$\S_n$, that moves $t_1$, and taking the product all the obtained
elements outputs $\tau_{\N}$.

Show that $E_n^{(n-1)}\subseteq T$. Take $e\in E_n^{(n-1)}$. Suppose
$$
e=\bigl\{A_1\cup A_1', \dots, A_k\cup A_k'\bigr\},
$$
where $k=\rank(e)\leq n-1$ and $|A_1|\geq 2$. Let $A_i=\{t_1^i,\dots
t_{m_i}^i\}, 1\leq i\leq k$. Construct the blocks $B_1,\dots, B_k$
as follows: $B_1=\{t_1^1\}$, $B_2$ consists of $|A_2|$ elements of
\begin{equation}\label{eq:is_ip} t_1^1, \dots, t_{m_1}^1, \dots,
t_1^k,\dots, t_{m_k}^k
\end{equation}
which follow $t_1^1$, $B_3$ consists of $|A_3|$ elements
of~\eqref{eq:is_ip} which follow the last element of $B_2$, and so
on, finally $B_k$ consists of the remaining $|A_k|+|A_1|-1$ elements
of~\eqref{eq:is_ip}. Set
$$
a=\{A_1\cup B_1', \dots, A_k\cup B_k'\}.
$$
The construction implies that some powers of $a$ and of $a^{-1}$ equal $\tau_{\N}$. Hence, $a,a^{-1}\in T$, and
thus $e=aa^{-1}\in T$.

Finally, since some power of every element of $\IP_n\setminus \S_n$
is an idempotent of $E_n^{(n-1)}\subset T$ and $T$ is isolated, we
have $\IP_n\setminus \S_n\subseteq T$. The statement follows.
\end{proof}

\subsection{Isolated subsemigroups of $\wPIP_n$}

\begin{theorem}\label{th:is_wpip}
The semigroups $\wPIP_n$, $\S_n$, $\wPIP_n\setminus \S_n$ and
$G(e)$, $e$ is an idempotent with $\corank(e)\leq 1$, and only them,
are isolated subsemigroups of $\wPIP_n$.
\end{theorem}

For the proof of this theorem we need some preparation.
The observation below follows from the definition of $\circ$.

\begin{lemma}\label{lem:rank_wpip}
Let $a\in\wPIP_n$. Then every block of $\dom(a^k)$ coincides with
some block of $\dom(a)$ and every block of $\ran(a^k)$ coincides
with some block of $\ran(a)$ for each $k\geq 1$.
\end{lemma}

Let $e\in E(\PIP_n)$.  Set $\corank(e)=|\codom(e)|=|\coran(e)|$.

\begin{lemma}\label{pr:is_wpip}
Let $e\in E(\wPIP_n)$ be such that $\corank(e)\leq 1$. Then $G(e)$
is an isolated subsemigroup of $\wPIP_n$.
\end{lemma}

\begin{proof}
Similarly to as in the proof of Proposition~\ref{pr:is_ip} it is
enough to prove that $a\in G(e)$ under the assumption that $a^k=e$
for some $k\geq 1$. Consider two possible cases.

{\em Case 1.} $\corank(e)=0$. Since $\coran(e)\supseteq \coran(a)$
and $\codom(e)\supseteq \codom(a)$ it follows that
$|\coran(a)|=|\codom(a)|=0$. Thus $\dom(a)$, $\dom(e)$, $\ran(a)$,
$\ran(e)$ are some partitions of $\N$. This and
Lemma~\ref{lem:rank_wpip} imply $\dom(a)=\dom(e)$ and
$\ran(a)=\ran(e)$. Therefore, $a\GH e$, implying $a\in G(e)$.

{\em Case 2.} $\corank(e)=1$. Assume that $\codom(e)$=$\{t\}$. By Lemma~\ref{lem:rank_wpip} there are
two possibilities: either $\dom(a)=\dom(e)$ and
$\ran(a)=\ran(e)$, or $\dom(a)=\dom(e)\cup\{t\}$ and
$\ran(a)=\ran(e)\cup\{t'\}$. In the first case we have $a\GH e$,
which yields $a\in G(e)$, as required. In the second case we would
have $a\GH f$ and then $e\in G(f)$, where $f$ is an idempotent
such that each generalised line of $e$ is a generalised line of $f$
and, besides, $f$ has the block $\{t,t'\}$, which is impossible.
\end{proof}

To proceed, we need to recall the description of isolated
subsemigroups of $\IS_n$ which is taken from~\cite[Chapter 5]{GM1}:

\begin{lemma}\label{lem:gan_k_is}
The semigroups $\IS_n$, $\S_n$, $\IS_n\setminus \S_n$, and $G(e)$,
where $e$ is an idempotent of rank $n-1$, and only them are isolated
subsemigroups of $\IS_n$.
\end{lemma}

\begin{proof}[Proof of Theorem~\ref{th:is_wpip}]
Applying Lemma~\ref{pr:is_wpip} and Theorem~\ref{th:wPIP-comp-isol}, it is
enough to prove the sufficiency. Let $T$ be an isolated subsemigroup
of $\wPIP_n$, such that $T\neq \S_n$ and $T\neq G(e)$ for any
idempotent $e$ of corank $0$ or $1$. We are to show that $T\supseteq
\wPIP_n\setminus \S_n$.

First show that $T$ has an idempotent of corank
at least $2$. Assume the converse. Then $T$ contains at least two
distinct idempotents $e$, $f$ such that $\corank(e)\leq 1$,
$\corank(f)\leq 1$. Since $ef\in T$ and $\corank (ef)\leq 1$, one of
$e,f$ must be equal to $ef$. Hence we can assume that $e\geq f$. We
have $G(e),G(f)\subseteq T$ by Lemma~\ref{pr:isol}. Observe that
among all the products of elements of $G(e)$ and $G(f)$ there are
elements  some powers of which are idempotents of corank at
least $2$.

Let $f\in T$ be an idempotent of corank at least $2$. Fix
$t_1,t_2\in\N$, $t_1\neq t_2$, such that $t_1, t_2\in \codom (f)$.
Define $a\in\wPIP_n$ as follows.
Each generalised line of $f$ is a generalised line of $a$. Besides,
$a$ has one more generalised line: $\{t_1,t_2'\}$. Then
$a^2=(a^{-1})^2=f$, the element  $\widetilde{f}=aa^{-1}$ is an
idempotent, and each generalised line of $f$ is a generalised line
of $\widetilde{f}$. In addition , $\widetilde{f}$ has exactly one
more line: $\{t_1,t_1'\}$. Since $T$ is isolated,
$G(\widetilde{f})\subseteq T$. Multiplying all the products of
elements from $G(\widetilde{f})$ by $f$ we obtain $0$. This shows
that $0\in T$.

Since $T\cap\IS_n\neq \varnothing$, it follows that $T\cap\IS_n$ is
an isolated subsemigroup of $\IS_n$, which by
Lemma~\ref{lem:gan_k_is} and $0\in T$ implies $\IS_n\setminus \S_n
\subseteq T$. Thus  $\wPIP_n\setminus\Sym_n\subseteq T$ by
Corollary~\ref{cor:new}.
\end{proof}

\subsection{Isolated subsemigroups of $\PIP_n$}

Let $Y\subset \N$ and $a\in \PIP_n$. We will call the set $Y$ {\em
invariant} with respect to $a$ if either $A\subset Y\cup Y'$ or
$A\cap (Y\cup Y')=\varnothing$ for each block $A$ of $a$. If $Y$ is
invariant with respect to $a$ denote by $a|_Y$ the element of
$\PIP_Y$ whose blocks are all blocks of $a$ which are contained in $Y\cup Y'$. The element $a|_Y$ will
be called the {\em restriction} of $a$ to $Y$. The semigroup $\IP_Y$
embeds into $\IP_n$ via the map sending $a\in\IP_Y$ to
the element of $\IP_n$ whose generalised lines are precisely the
generalised lines of $a$, and all the other blocks are points. We
will identify $\IP_Y$ with its image under this embedding.

\begin{lemma}\label{pr:is_pip}
Let $n\geq 3$. The semigroups \begin{enumerate}[1)]
\item $\IP_n$,
$\IP_n\setminus \S_n$, $\S_n$, $G(e)$, where $e$ is an idempotent of
rank $n-1$ of $\IP_n$;
\item $\IP_Y$, $\IP_Y\setminus \S_Y$, $\S_Y$, $G(e)$, where $e$ is
an idempotent of rank $n-2$ of $\IP_Y$, where
$Y=\N\setminus\{t\}$, $t\in\N$;
\item $\PIP_n$, $\PIP_n\setminus \S_n$
\end{enumerate}
are isolated subsemigroups of $\PIP_n$.
\end{lemma}
\begin{proof}
The proof is a straightforward verification. It resembles the proofs
of Proposition~\ref{pr:is_ip} and Lemma~\ref{pr:is_wpip}.
\end{proof}

\begin{theorem}\label{th:is_pip}
Let $n\geq 3$. The semigroups listed in Lemma~\ref{pr:is_pip} and
only them are isolated subsemigroups of $\PIP_n$.
\end{theorem}

\begin{proof}
Let $T$ be an isolated subsemigroup of $\PIP_n$. If $T\subset \IP_n$
then $T$ must be an isolated subsemigroup of $\IP_n$. Therefore,
applying Theorem~\ref{th:is_ip}, we see that $T$ is one of the
semigroups listed in the first item of Lemma~\ref{pr:is_pip}.

Suppose $T\setminus \IP_n\neq \varnothing$. Then $T$ contains an
idempotent of corank $1$ (this can be shown using arguments similar
to those from the third paragraph of the proof of
Theorem~\ref{th:is_wpip}, where an idempotent $\widetilde{f}$ is
being constructed by $f$). It follows that there is $Y\subset \N$,
$Y=\N\setminus\{t\}$, $t\in \N$, such that $T\cap \IP_Y
\neq\varnothing$. It follows that $T\cap \IP_Y$ is an isolated
subsemigroup of $\IP_Y$. If $T\subseteq \IP_Y$ then $T$ is one of
the semigroups of the second item of Lemma~\ref{pr:is_pip}.

Suppose that $T\setminus \IP_Y\neq\varnothing$. Then $T$ has at
least two idempotents $e$ and $f$ such that there is no proper
subset $Z$ of $\N$ for which $e,f\in\IP_Z$. Since $e,f,ef\in T$ it follows that we may assume $e> f$. Now, $G(e), G(f)\subset
T$ imply $0\in T$. Hence $T\cap \IS_n$ is an isolated subsemigroup
of $\IS_n$ containing the zero. This and Lemma~\ref{lem:gan_k_is}
show that $\IS_n\setminus \S_n\subseteq T$.

To complete the proof show that $\PIP_n\setminus
\S_n\subseteq T$. It is enough to show that
$\widetilde{E}_n^{(n-1)}\subseteq T$. Let $e\in \PIP_n\setminus
\IS_n$ be an idempotent. Let  $Z=\N\setminus\codom(e)$. If $Z=\N$ then $e\in T$ by arguments at
the end of the proof of Theorem~\ref{th:is_ip}. Let $\N\setminus Z\neq\varnothing$. We have that
$e|_Z\in E(\IP_Z\setminus \S_Z)$. We claim that it is enough to show that the element $\widetilde{\tau}_Z$, having
the only generalised line $Z\cup Z'$ and all the other blocks
points, belongs to $T$. Indeed, if $\widetilde{\tau}_Z\in T$ then applying the arguments similar to those at
the end of the proof of Theorem~\ref{th:is_ip}, we obtain that
$e|_Z\in T|_Z$, implying that $f\in T$ for some $f\in\PIP_n$ with $e|_Z=f|_Z$. Since we also know  that $1|_Z\in \IS_n\setminus\S_n\subseteq T$, we have that $e=1|_Z f\in T$ as well. Take $t\in Z$. Set $a$ to be the
element of $\PIP_n$ with the only one generalised line $Z\cup
\{t'\}$, and all the other blocks points. Then $a^2=(a^{-1})^2=0$,
while $aa^{-1}=\widetilde{\tau}_Z$. The statement follows.
\end{proof}

\section{Automorphisms of $\PIP_X$ and $\wPIP_X$}\label{sec:aut}

\subsection{Automorphisms of $\PIP_X$}

Let $Y\subset X$. We will need to consider the following
subsemigroups of $\PIP_X$:
\begin{equation*}
\oS_Y=\bigl\{a\in\S_X:~a~\mbox{contains the blocks}~\{t,t'\},~t\in
X\setminus Y\bigr\},
\end{equation*}
\begin{equation*}
\oIS_Y=\bigl\{a\in\IS_X:~a~\mbox{contains the blocks}~\{t,t'\},~t\in
X\setminus Y\bigr\}~\mbox{and}
\end{equation*}
\begin{equation*}
\oIP_Y=\bigl\{a\in\IP_X:~a~\mbox{contains the blocks}~\{t,t'\},~t\in
X\setminus Y\bigr\}.
\end{equation*}

Let $\Aut(S)$ denote the group of automorphisms of a semigroup $S$.

\begin{theorem}\label{th:Aut-PIP}
$\Aut(\PIP_X)\cong\Sym_X$. Moreover, for every
$\varphi\in\Aut(\PIP_X)$ there is $\pi\in \S_X$ such that
$a^{\varphi}=\pi^{-1}a\pi$, $a\in \PIP_X$.
\end{theorem}

\begin{proof}
Let $\varphi\in\Aut(\PIP_X)$. Take $x\in X$. Since $\Sym_X$ is the
group of units of $\PIP_X$, in should be preserved by $\varphi$:
$\varphi(\Sym_X)=\Sym_X$. For $u\in\PIP_X$ and a subsemigroup
$T\subseteq \PIP_X$ let
$$
\St_{T}^r(u)=\{s\in T\mid us=u\}, \,\, \St_{T}^l(u)=\{s\in T\mid
su=u\}.
$$
Recall that for $x\in X$ by $\alpha_{x}$ we denote the idempotent
$\bigl\{\{t,t'\}_{t\in X\setminus\{x\}}\bigr\}$
$\in\PIP_X$.

Observe that for an idempotent $u\in\PIP_X$ $|\St_{\S_X}^r(u)|=1$ if
and only if $u=\alpha_z$ for some $z\in X$. It follows that for each
$x\in X$ there is $g(x)\in X$ such that
$\varphi(\alpha_x)=\alpha_{g(x)}$. This defines a permutation $g\in
\Sym_X$.

Show that $\varphi(\IS_X)=\IS_X$. Let $a=\bigl\{\{t,\pi(t)'\}_{t\in
I}\bigr\}\in\IS_X$, where $\pi:~I\to\pi(I)$ is a
bijection. For all $z\in X\setminus I$ and $r\in X\setminus \pi(I)$
we have $\alpha_z a=a=a\alpha_r$. Passing in this equality to
$\varphi$-images, we see that $\varphi(a)$ should contain the blocks
$\{q\}$, $q\in g(X\setminus I)$, and $\{r'\}$, $r\in
g(X\setminus\pi(I))$. Let $t_0\in I$. Notice that the equality
\begin{equation}\label{eq:aaa}
\alpha_{t_0}a=\alpha_{z}\cdot \alpha_{t_0}a\cdot \alpha_{r}
\end{equation}
holds if and only if $z\in (X\setminus
I)\cup\{t_0\}~\mbox{and}~r\in\bigl(X\setminus\pi(I)\bigr)\cup\{\pi(t_0)\}$.
Going in~\eqref{eq:aaa} to $\varphi$-images, we obtain 
$$\alpha_{g(t_0)}\varphi(a)=\alpha_{g(z)}\cdot
\alpha_{g(t_0)}\varphi(a)\cdot \alpha_{g(r)}.$$
Similarly as above we have that the equality
$$\alpha_{g(t_0)}\varphi(a)=\alpha_{z}\cdot
\alpha_{g(t_0)}\varphi(a)\cdot \alpha_{r}$$
holds if and only if $z\in
g\bigl((X\setminus I)\cup\{t_0\}\bigr)$ and $r\in
g\bigl(\bigl(X\setminus\pi(I)\bigr)\cup\{\pi(t_0)\}\bigr)$. The
latter implies that $\varphi(a)$ contains a block
$\{g(t_0),g(\pi(t_0))\}$. Now we can assert that
$\varphi(a)=\bigl\{\{g(t),g(\pi(t))'\}\bigr\}_{t\in I}$. It
follows that $\varphi(\IS_X)=\IS_X$. Moreover, for every $Y\subset
X$ we have
\begin{equation}\label{eq:IS}
\varphi(\oIS_Y)=\oIS_{g(Y)}.
\end{equation}

Show that $\varphi(\IP_X)=\IP_X$. Observe that the elements of
$\IP_X$ may be characterized as follows: $b\in\IP_X$ if and only if
$\alpha_x b\ne b$ and $b\alpha_x \ne b$ for all $x\in X$. Let
$b=\bigl\{(A_i\cup B_i')_{i\in I}\bigr\}\in\IP_X$. The equality
$\alpha_u b=\alpha_v b$ holds if and only if $u$ and $v$ belong to
$A_i$ for some $i\in I$, the equality $b \alpha_u =b \alpha_v$ holds
if and only if $u$ and $v$ belong to $B_i$ for some $i\in I$, and
the equality $\alpha_u b= b\alpha_v$ holds if and only if $u\in A_i$
and $v\in B_i$ for some $i\in I$. Going to $\varphi$-images and
using the fact that $\varphi(\alpha_x)=\alpha_{g(x)}$, $x\in X$, we
can assert that $\varphi(b)=\bigl\{\bigl(g(A_i)\cup
g(B_i)'\bigr)_{i\in I}\bigr\}$. Thus $\varphi(\IP_X)=\IP_X$ and,
moreover,
\begin{equation}\label{eq:IP}
\varphi(\oIP_Y)=\oIP_{g(Y)}
\end{equation}
for every $Y\subset X$. Since $\oS_Y=\oIS_Y\cap \oIP_Y$,
applying~\eqref{eq:IS} and~\eqref{eq:IP} we obtain
\begin{equation}\label{eq:S}
\varphi(\oS_Y)=\varphi(\oIS_Y)\cap \varphi(\oIP_Y)=\oS_{g(Y)}.
\end{equation}

Let $a=\bigl\{(U_i\cup V_i')_{i\in
I}\bigr\}\in\PIP_X$.
Observe that
\begin{equation*}
\St_{\IP_X}^l(a)= (\oIP_{X\setminus\bigcup\limits_{i\in
I}U_i})\oplus(\bigoplus_{i\in I}\oIP_{U_i}); \,\,
\St_{\IP_X}^r(a)=(\oIP_{X\setminus\bigcup\limits_{i\in
I}V_i})\oplus(\bigoplus\limits_{i\in I}\oIP_{V_i});
\end{equation*}
\begin{equation*}
\St_{\IS_X}^l(a)=(\oIS_{X\setminus\bigcup\limits_{i\in I}U_i})
\oplus(\bigoplus\limits_{i\in I}\oS_{U_i});\,\,
\St_{\IS_X}^r(a)=(\oIS_{X\setminus\bigcup\limits_{i\in I}V_i})
\oplus(\bigoplus\limits_{i\in I}\oS_{V_i}).
\end{equation*}
We observe that the equalities
$$\St_{\IP_X}^l(a)=\St_{\IP_X}^l(b),
\St_{\IP_X}^r(a)=\St_{\IP_X}^r(b\,),
\St_{\IS_X}^l(a)=\St_{\IS_X}^l(b),
\St_{\IS_X}^r(a)=\St_{\IS_X}^r(b)$$ hold for some $b\in \PIP_X$ if
and only if $\dom(a)=\dom(b)$ and $\ran(a)=\ran(b)$, which by
Proposition~\ref{pr:Green_and_subgroups}, is
equivalent to $a\GH b$.

By~\eqref{eq:IS},~\eqref{eq:IP} and~\eqref{eq:S} we have
\begin{equation}\label{eq:aut}
\varphi(a)=\bigl\{g(U_i)\cup g(V_{\pi(i)})'\bigr\}_{i\in
I}
\end{equation}
for some bijection $\pi:~I\to I$. Let us show that $\pi$ should be
the identity map. Let $j\in I$. Fix $u_j\in U_j$. We compute
$$
\alpha_{u_j}a=\bigl\{(U_i\cup V_i')\bigr\}_{i\in
I\setminus\{j\}}.
$$
By~\eqref{eq:IS},~\eqref{eq:IP} and~\eqref{eq:S} we have
\begin{equation}\label{eq:final}
\coran(\varphi(\alpha_{u_j}a))=\bigl\{\{t'\},
t\not\in\bigcup\limits_{i\in I}V_i, \{t'\}, t\in g(V_j)\bigr\}.
\end{equation}
From the other hand,
$\varphi(\alpha_{u_j}a)=\alpha_{g(u_j)}\varphi(a)$, and thus
\begin{equation}\label{eq:!final}
\coran(\varphi(\alpha_{u_j}a))=\bigl\{\{t'\},
t\not\in\bigcup\limits_{i\in I}V_i, \{t'\}, t\in
g(V_{\pi(j)})\bigr\}.
\end{equation}
It follows from~\eqref{eq:final} and~\eqref{eq:!final} that
$\pi(j)=j$, and then $\pi$ is the identity map. Hence
$\varphi(a)=g^{-1}ag$, $a\in\PIP_X$. The proof is completed.
\end{proof}

\subsection{Automorphisms of $\wPIP_X$}

Let $Y\subseteq X$. Set $\varepsilon_Y=\bigl\{Y\cup Y'\bigr\}.$ The element
$\varepsilon_Y$ is an idempotent of rank $1$. If $\varepsilon$ is an
idempotent of rank $1$, denote by $Y(\varepsilon)$ such a subset
$Y\subseteq X$ that $\varepsilon_{Y(\varepsilon)}=\varepsilon$.

\begin{theorem}\label{th:Aut-wPIP}
$\Aut(\wPIP_X)\cong\Sym_X$.
\end{theorem}

\begin{proof}
Let $\varphi\in\Aut(\wPIP_X)$. The maps $\varepsilon_Y \mapsto Y$
and $Y\mapsto \varepsilon(Y)$ are mutually inverse bijections
between the idempotents of rank $1$ of $\wPIP_X$ and nonempty
subsets of $X$. It follows that $\varphi$ induces some permutation
$\pi$ on $2^X\setminus \{\varnothing\}$.

Show that $A\cap B=\varnothing$ implies $\pi(A)\cap
\pi(B)=\varnothing$ for all $A,B\subseteq X$. Consider the
idempotent $e=\bigl\{A\cup A', B\cup B'\bigr\}$. Let
$f=\varphi(e)=\{C\cup C', D\cup D'\}$ ($\rank(f)=2$ because
$\rank(e)=2$, and ranks of  idempotents are preserved by
automorphisms as they may be characterised in terms of the natural
order). Since $\varepsilon_Ae=\varepsilon_A$ and
$\varepsilon_Be=\varepsilon_B$, going to $\varphi$-images, we obtain
$\varepsilon_{\pi(A)}f=\varepsilon_{\pi(A)}$ and
$\varepsilon_{\pi(B)}f=\varepsilon_{\pi(B)}$. It follows that $f$
has the blocks $\pi(A)\cup \pi(A)'$ and $\pi(B)\cup \pi(B)'$. Taking
into account that $\rank(f)=2$, we see that
$\{C,D\}=\{\pi(A),\pi(B)\}$. Since $C\cap D=\varnothing$, than also
$\pi(A)\cap \pi(B)=\varnothing.$

Show now that $\pi$ maps one-element subsets of $X$ to one-element
subsets. Assume the converse. Let $x\in X$ be such that
$\pi(\{x\})=M$, where $|M|\geq 2$. Take $y, z\in M$, $y\neq z$. Let
$M_y$ and $M_z$ denote the sets satisfying $\pi(M_y)=\{y\}$ and
$\pi(M_z)=\{z\}$, respectively. Since $\{y\}\cap\{z\}=\varnothing$,
by the argument from the previous paragraph we obtain $M_y\cap
M_z=\varnothing$. On the other hand, using $\{y\}\cap M\neq
\varnothing$ and $\{z\}\cap M\neq \varnothing$, we obtain that it
must be $M_y\cap \{x\}\neq\varnothing$ and $M_z\cap
\{x\}\neq\varnothing$. But then $x\in M_y\cap M_z$, which is
impossible. The restriction of $\pi$ to one-element subsets of $X$
defines a permutation $g\in \S_X$.

We proceed by showing that $\pi(M)=g(M)=\{g(m)\mid m\in M\}$ for
each subset $M$ of $X$. Indeed, since $M\cap \{t\}=\varnothing$,
$t\in X\setminus M$, it follows that $\pi(M)\subseteq g(M)$. Similar
arguments applied for the automorphism $\varphi^{-1}$ ensure that
$\pi^{-1}(g(M))\subseteq M$, and thus $g(M)\subseteq \pi(M)$. The
reverse inclusion is established similarly.

Let $a\in \wPIP_X$. Suppose that $a$ has a block $A\cup B'$. Show
that $\varphi(a)$ has the block $g(A)\cup g(B)'$. Indeed,
$\varepsilon_A a\varepsilon_B\neq 0$. Going to $\varphi$-images, we
obtain $\varepsilon_{g(A)} \varphi(a) \varepsilon_{g(B)}\neq 0$. The
latter implies that $\varphi(a)$ has the block $g(A)\cup g(B)'$, as
required. It follows that $A\cup B'$ is a generalised line
of $a$ if and only if $g(A)\cup g(B)'$ is a generalised line of
$\varphi(a)$, which completes the proof.
\end{proof}

\section{$\PIP_n$ and $\wPIP_n$ are embeddable into $\IS_{2^n-1}$}\label{sec:ef}

Let $S$ be an inverse semigroup with the natural partial order
$\varrho$ on it. The following definitions are taken from
\cite[p.~188]{Howie}. An inverse subsemigroup $H$ of $S$ is called a
\emph {closed inverse subsemigroup} of $S$ if $H\varrho=H$. Let
\begin{equation}
\mathcal{C}=\mathcal{C}_H=\bigl\{(Hs)\varrho:~ss^{-1}\in H\bigr\}
\end{equation} be the set of all \emph{right $\varrho$-cosets} of $H$.

Let, further, \begin{equation}
\phi_H(s)=\bigr\{\bigl((Hx)\varrho,(Hxs)\varrho\bigr):~
(Hx)\varrho,(Hxs)\varrho\in\mathcal{C}\bigl\}
\end{equation}
be the \emph{effective transitive representation}
$\phi_H:S\to\IS_{\mathcal{C}}$. If $K$ and $H$ are two closed
inverse subsemigroups of $S$, the representations $\phi_K$ and
$\phi_H$ are \emph{equivalent} if and only if there exists $a\in S$
such that $a^{-1}Ha\subseteq K$ and $aKa^{-1}\subseteq H$ (see
\cite[Proposition IV.4.13]{Petrich}).

\begin{theorem}\label{th:ef}
Let $n\geq 2$. Up to equivalence, there is only one faithful
effective transitive representation of $\PIP_n$ (respectively
$\wPIP_n$), namely to $\IS_{2^n-1}$. In particular, $\PIP_n$ and
$\wPIP_n$ embed into $\IS_{2^n-1}$.
\end{theorem}

\begin{proof}
We prove the statement for the case of $\PIP_n$, the other case
being treated analogously. Suppose $H$ is a closed inverse
subsemigroup of $\PIP_n$. Denote by $\omega$ the natural partial
order on $\PIP_n$. First we observe that $H=G\omega$ for some
subgroup $G$ of $\PIP_n$. Indeed, since $\PIP_n$ is finite,  $E(H)$
contains a zero element. It remains to apply \cite[Proposition
IV.5.5]{Petrich}, which claims that if the set of idempotents of a
closed inverse subsemigroup contains a zero element, then this
subsemigroup is a closure of some subgroup of the original
semigroup. Denote by $f$ the identity element of the group $G$.

Now we prove that if $f=0$ then $\phi_H$ is not faithful. We have
$H=0\omega=\PIP_n$ and hence $(Hx)\omega\supseteq 0\omega=\PIP_n$
for all $x\in\PIP_n$. Thus $(Hx)\omega=\PIP_n$ for all
$x\in\PIP_n$. Then $\lmod\phi_H(\PIP_n)\rmod=1$ and so $\phi_H$ is
not faithful.

Let now $\rank(f)\geq 2$. We will show that in this case $\phi_H$ is
not faithful either. Take $b\in D_1$ where $D_1$ denotes the set of
elements of $\PIP_n$ of rank $1$. Since $bb^{-1}\in D_1$ we have
that $bb^{-1}\notin H$ and therefore $(Hb)\omega\notin\mathcal{C}$.
This implies that $\phi_H(b)$ is equal to the zero element of
$\IS_{\mathcal{C}}$. Then due to $\lmod D_1\rmod\geq 2$ we obtain
that $\phi_H$ is not faithful.

Let finally $\rank(f)=1$. We will show that in this case $\phi_H$ is
faithful. Observe that $H=f\omega$. Let
$f=\varepsilon_E=\bigl\{E\cup E'\bigr\}$ where $E\neq\varnothing$.
Suppose that $\phi_H(s)=\phi_H(t)$ for some $s$ and $t$ from
$\PIP_n$. Without loss of generality assume that $s\neq 0$. Suppose
that $s$ contains a block $A\cup B'$. Consider the element
$x=\bigl\{E\cup A'\bigr\}$.
Then $(Hx)\omega$ and $(Hxs)\omega$ belong to $\mathcal{C}$. This
implies that $(Hxs)\omega=(Hxt)\omega$. The latter means that $t$
contains some generalised lines whose union is the block $A\cup B'$.
Changing the roles of $s$ and $t$ we obtain that both $s$ and $t$
contain the block $A\cup B'$. Thus $s=t$, as required.

Observe that all the idempotents of $\PIP_X$ of rank $1$ are
precisely the primitive idempotents. Let $g$ be a primitive
idempotent of $\PIP_n$. We will show that
$\lmod\mathcal{C}_{g\omega}\rmod=2^n-1$. Note that
$\mathcal{C}_{g\omega}=\bigl\{(gs)\omega:~ss^{-1}\geq g\bigr\}$.
We have $(gs)\omega=(gt)\omega$ if and only if $gs=gt$, that is, the
number of different sets $(gs)\omega$, $ss^{-1}\geq g$, is equal to
the number of different nonempty subsets of $\N$, which equals
$2^n-1$.

To complete the proof we note that for two primitive idempotents
$f_1,f_2\in\PIP_n$ we have that $\phi_{f_1\omega}$ and
$\phi_{f_2\omega}$ are equivalent by the definition of equivalent
representations.
\end{proof}

\noindent G.K.: Algebra, Department of Mathematics and Mechanics,
Kyiv Taras
Shevchenko University, 64 Volodymyrska st., 01033 Kyiv, UKRAINE,\\
e-mail: {\tt akudr\symbol{64}univ.kiev.ua} and

\noindent Centre for Systems and Information Technologies, University of Nova Gorica,
Vipavska 13, PO Box 301, Rozna Dolina, SI-5000 Nova Gorica, SLOVENIA,\\ 
e-mail: {\tt ganna.kudryavtseva\symbol{64}p-ng.si} 
\vspace{0.3cm}

\noindent V.M.: School of Mathematics and Statistics, University
of St Andrews, St Andrews, Fife, KY 16 9SS, SCOTLAND,\\
e-mail: {\tt victor\symbol{64}mcs.st-and.ac.uk}

\end{document}